\documentclass{amsart}
\usepackage{amssymb,euscript,amsmath, mathrsfs}
\usepackage[dvips]{graphicx}
\usepackage[dvips]{color}

\newcounter{ENUM}
\newcommand{\itm}{\item}
\newenvironment{ilist}{\renewcommand{\theENUM}{\roman{ENUM}}\renewcommand{\itm}{\addtocounter{ENUM}{1}\item[(\theENUM)]}\begin{itemize}\setcounter{ENUM}{0}}{\end{itemize}}

\newenvironment{nlist}[1][0]{\renewcommand{\theENUM}{\arabic{ENUM}}\renewcommand{\itm}{\addtocounter{ENUM}{1}\item[(\theENUM)]}\begin{itemize}\setcounter{ENUM}{#1}}{\end{itemize}}

\newcommand{\margh}[1]{}

\def\risom{\overset{\sim}{\rightarrow}}

\input xy
\xyoption{all}
\CompileMatrices

\def\N{{\mathbb N}}
\def\A{{\mathbb A}}
\def\P{{\mathbb P}}

\def\fF{{\mathbb F}}
\def\F{{\mathscr F}}

\def\L{{\mathscr L}}
\def\M{{\mathscr M}}
\def\E{{\mathscr E}}
\def\O{{\mathscr O}}

\def\cHom{\mathcal{H}om}
\def\cEnd{\mathcal{E}nd}

\def\i{{\mathfrak i}}

\def\m{{\mathfrak m}}

\def\vp{\varphi}

\def\Nc{\nabla ^{\text{can}}}
\def\Id{\operatorname{I}}

\def\Ext{\operatorname{Ext}}

\def\Hom{\operatorname{Hom}}

\def\End{\operatorname{End}}
\def\Aut{\operatorname{Aut}}

\def\Tr{\operatorname{Tr}}
\def\Spec{\operatorname{Spec}}
\def\Pic{\operatorname{Pic}}

\def\Res{\operatorname{Res}}
\def\Conn{\operatorname{Conn}}

\def\ind{\operatorname{ind}}

\def\ns{\operatorname{ns}}

\def\nod{\operatorname{nod}}

\numberwithin{equation}{section}
\newtheorem{thm}{Theorem}[section]
\newtheorem{prop}[thm]{Proposition}
\newtheorem{lem}[thm]{Lemma}
\newtheorem{cor}[thm]{Corollary}

\theoremstyle{definition}
\newtheorem{defn}[thm]{Definition}
\newtheorem{ques}[thm]{Question}

\newtheorem{sit}[thm]{Situation}
\theoremstyle{remark}
\newtheorem{notn}[thm]{Notation}

\begin{document}
\title{{F}robenius-Unstable Bundles and $p$-Curvature}
\author{Brian Osserman}
\begin{abstract}
We use the theory of $p$-curvature of connections to analyze stable vector
bundles of rank 2 on curves of genus 2 which pull back to unstable bundles
under the Frobenius morphism. We take two approaches, first using explicit
formulas for $p$-curvature to analyze low-characteristic cases, and then
using degeneration techniques to obtain an answer for a general curve by
degenerating to an irreducible rational nodal curve, and applying the
results of \cite{os10} and \cite{os7}. We also apply our explicit formulas
to give a new description of the strata of curves of genus 2 of different 
$p$-rank.
\end{abstract}
\thanks{This paper was partially supported by fellowships from the National
Science Foundation and Japan Society for the Promotion of Science.}
\maketitle

\section{Introduction}

The primary theme of this paper is to use the following question as an invitation to a detailed study of the theory of $p$-curvature of connections:

\begin{ques} Given a smooth curve $C$ of genus 2 over an algebraically
closed field $k$ of positive characteristic, what is the number of {\bf
Frobenius-unstable} vector bundles of rank 2 and trivial determinant on
$C^{(p)}$? That is, if $F:C \rightarrow C^{(p)}$ denotes the relative
Frobenius morphism from $C$ to its $p$-twist, how many vector bundles $\F$ are there on $C^{(p)}$ (of the stated rank and determinant) which are themselves semistable, but for which $F^* \F$ is unstable?
\end{ques}

Because semistability is preserved by pullback under separable morphisms
(see \cite[Lem. 3.2.2]{h-l}), the Frobenius-unstable case is in some sense a universal case for destabilization. Furthermore, Frobenius-unstable bundles are closely related to the study of the generalized Verschiebung, and its relationship to $p$-adic representations of the fundamental group of $C$, in the case that $C$ is defined over a finite field; see \cite{os9} for details.

The analysis of our question is in two parts: first, we use explicit formulas for $p$-curvature to calculate the answer directly for odd characteristics $\leq 7$; and second, we use the abstract theory of $p$-curvature to give a new proof of the answer for a general curve of genus 2 in any odd characteristic, via degeneration to an irreducible rational nodal curve and application of the results of \cite{os10} and \cite{os7}. The latter result is originally due to Mochizuki; see \cite{mo3} and 
\cite{os6}. The main advantage of the explicit approach, as compared to the more general degeneration argument, is that the $p$-curvature 
formulas may be used to study arbitrary smooth curves, and do not give 
results only for general curves. This distinction is underscored by an 
algorithm derived via the same techniques to explicitly describe the loci 
of curves of genus 
$2$ and $p$-ranks $0$ or $1$ in any specified characteristic. Additionally, the explicit approach is useful for computing examples in order to formulate conjectures; one aim of this paper is therefore to serve as an illustration of how $p$-curvature may be used very concretely for experimental purposes, and more theoretically for more general results.

Our main theorem is:

\begin{thm}\label{exp-main} Let $C$ be a smooth, proper curve of genus $2$ 
over an algebraically closed field $k$ of characteristic $p$; it may be
described on an affine part by $y^2=g(x)$ for some quintic $g$. Then
the number of semistable vector bundles on $C$ with trivial determinant 
which pull
back to unstable vector bundles under the relative Frobenius morphism is:
\begin{itemize}
\item[$p=3$:] $16 \cdot 1$;
\item[$p=5$:] $16 \cdot e_5$, where $e_5=5$ for $C$ general, and is given
for an arbitrary $C$ as the number of distinct roots of a
quintic polynomial with coefficients in terms of the coefficients of $g$;
\item[$p=7$:] $16 \cdot e_7$, where $e_7=14$ for $C$ general, and is 
given for an arbitrary $C$ as the number of points in the intersection of 
four curves in $\A^2$ whose coefficients are in terms of the coefficients of 
$g$.
\item[$p>2$:] (Mochizuki \cite{mo3}, \cite{os6}) $16 \cdot \frac{p^3-p}{24}$ for $C$ general.
\end{itemize}

Furthermore, when $C$ is general, any Frobenius-unstable bundle $\F$ has no non-trivial deformations which yield the trivial deformation of $F^* \F$.
\end{thm}

There is a considerable amount of literature on Frobenius-unstable vector bundles. Gieseker and Raynaud produced certain examples of
Frobenius-unstable bundles in \cite{gi2} and \cite[p. 119]{ra2}, but, aside from the results of Mochizuki discussed below, the first classification-type result is due to Laszlo and Pauly, who answered our main question in characteristic 2: there is always a single Frobenius-unstable bundle (see \cite{l-p}, argument for Prop. 6.1 2.; the equations for an ordinary
curve are not used). Joshi, Ramanan, Xia and Yu obtain results on the
Frobenius-unstable locus in characteristic $2$ for higher-genus curves in
\cite{j-r-x-y}.
Most recently, and concurrently with the initial preparation
of the present paper, Lange and Pauly \cite{l-p3} have recovered the formula of Theorem \ref{exp-main} for general $C$ in 
the case of ordinary curves via a completely different approach, although
they obtain only an inequality, rather than an equality.

However, the most comprehensive results to date follow from Mochizuki's work (see \cite{mo3} and \cite{os6}), which was carried out in the context of $\P^1$-bundles on curves in any odd characteristic, via degeneration techniques
quite similar to those which we pursue in Sections 8 and 9. Indeed, key results and their arguments in Sections \ref{s-exp-det}, \ref{s-def-background}, and \ref{s-def-deform} are essentially the same as Mochizuki's; in the first case, the argument presented here was discovered independently, while in the other cases, the author's original arguments were more complicated and less general than Mochizuki's, and have thus been replaced. There are several justifications for the logical redundancy: the arguments in question are all quite short, and it seems desirable to have a self-contained proof of the main theorem, without translating to projective bundles and back; the argument of Section \ref{s-exp-det} is actually substantially simpler in our case of curves of genus $2$; and finally, the gluing statements of Section \ref{s-def-background} require some ridigifying hypotheses in the context of vector bundles that do not arise in Mochizuki's work. 

Lastly, we remark that as discussed in 
\cite{os6}, Mochizuki's strategy is to degenerate to totally degenerate curves, while our strategy is to degenerate to irreducible nodal curves. Aside from allowing one to make more naive arguments in terms of explicit degenerations, ours is a substantially more difficult approach, since after reducing the problem to self-maps of $\P^1$ with prescribed ramification, in Mochizuki's case it suffices to handle the case of three ramification points, while our argument requires four, and is therefore far more complicated; see \cite{os7} for details. However, degenerating to irreducible curves is helpful for studying Frobenius-unstable bundles in higher genus; see
\cite{os6}.

We begin in Section \ref{s-exp-background} by relating our main question to $p$-curvature, and Section \ref{s-exp-gen-pcurve} is then devoted to
developing explicit and completely general combinatorial formulas for
$p$-curvature. We make certain necessary computations for genus $2$ curves in Section \ref{s-exp-ftheta}, which we also apply to obtain an explicit algorithm for generating $p$-rank formulas in any
given odd characteristic. Section \ref{s-exp-conns} is devoting to computing the space of connections on a certain unstable bundle, and in Section \ref{s-exp-pcurve} we conclude the 
computation with explicit descriptions of the locus of vanishing
$p$-curvature in characteristics $3, 5 \text{ and } 7$. The space of
connections on the same bundle having nilpotent $p$-curvature is shown to be finite and flat in Section \ref{s-exp-det}, again by explicit computation; this completes
the proof of Theorem \ref{exp-main} for $p\leq 7$, and also provides a key step of the general case. In Section \ref{s-def-background} we discuss the relationship between connections on nodal curves and their normalizations, and finally in Section \ref{s-def-deform} we show that connections on nodal curves deform, and apply the results of \cite{os10} and \cite{os7} to conclude our main theorem.

Computations were carried out in Maple and Mathematica, and in the case of the $p$-curvature formulas of Section \ref{s-exp-gen-pcurve}, using simple C code.

The contents of this paper form a portion of the author's 2004 PhD thesis 
at MIT, under the direction of Johan de Jong.

\section*{Acknowledgements}

I would like to thank Johan de Jong for his tireless and invaluable
guidance. I would also like to thank Shinichi Mochizuki, Ezra Miller, David Helm, and Brian Conrad for their helpful conversations. 

\section{From {F}robenius-instability to $p$-curvature}\label{s-exp-background}

We begin by explaining how classification of Frobenius-unstable vector
bundles is related to $p$-curvature of connections. For the basic theory of
connections and $p$-curvature, we refer the reader to \cite[\S 1, \S
5]{ka1}. Note that the induced connection on tensor products descends to
wedge products, so that for a vector bundle $\E$ with connection, we obtain
an induced {\bf determinant connection} on $\det \E$. Additionally, given
$\vp \in \Aut(\E)$ and a $\nabla$ on $\E$, we refer to the operation of
conjugation by $\vp$ on $\nabla$ as {\bf transport}. We summarize the basic
results relating Frobenius with $p$-curvature, due to to Katz \cite{ka1}.

\begin{thm}\label{exp-cartier-vect} Let $X$ be a smooth $S$-scheme, with 
$S$ having characteristic $p$, and let $F: X \rightarrow X^{(p)}$ be the 
relative Frobenius morphism. Then for any vector bundle $\F$ on $X^{(p)}$, $F^* \F$ is equipped with a canonical connection $\Nc$. For any vector bundle $\E$ with connection $\nabla$ on $X$, the kernel of $\nabla$, denoted $\E^\nabla$, is naturally an $\O_{X^{(p)}}$-module. 

The operations $\F \mapsto (F^* \F, \Nc)$ and $(\E, \nabla) \mapsto \E^{\nabla}$, are mutually inverse functors, giving an equivalence of categories between the category of vector bundles of rank $n$ on $X^{(p)}$ 
and the full subcategory of the category of vector bundles of rank $n$ with 
integrable connection on $X$ consisting of objects whose connection has 
$p$-curvature zero.

Furthermore, the same statement holds when restricted to the full 
subcategories of vector bundles with trivial determinant on $X^{(p)}$, and 
vector bundles with connection both having trivial determinant on $X$.
\end{thm}

\begin{proof} See \cite[\S 5]{ka1}, and in particular \cite[Thm. 5.1]{ka1}. It only remains to check that the categorical equivalence on coherent sheaves gives an equivalence on vector bundles, and again in the case of trivial determinant. The first assertion follows from the fact that $F$ is faithfully flat when $X/S$ is smooth. The second is easily checked by verifying that the operation $\F \mapsto (F^*\F, \Nc)$ commutes with taking determinants.
\end{proof}

Thus, $p$-curvature is naturally related to the study of
Frobenius-pullbacks. The categorical equivalence implies that isomorphism
classes of $\F$ will correspond to transport equivalence classes of
connections with vanishing $p$-curvature. In the case of our particular
question, the relationship is particularly helpful. We assume we are in the following situation.

\begin{sit}\label{exp-genus} $C$ is a smooth, proper curve of genus 2, 
over an algebraically closed field $k$ of characteristic $p$.
\end{sit}

In this situation, Joshi and Xia showed that there are at most finitely 
many Frobenius-unstable vector bundles of rank $2$ and trivial determinant 
on $C$ (see \cite[Thm. 3.2]{j-x}, although we will also 
obtain a more direct proof from Corollary \ref{exp-finite}), and also gave
the following description of them (see \cite[Prop. 3.3]{j-x}):

\begin{prop}\label{exp-unstable} (Joshi-Xia) Let $\F$ be a semistable rank $2$ 
vector bundle on $C$ with 
trivial determinant, and suppose $\E = F^* \F$ is unstable. Then there is a
non-split exact sequence 
$$0 \rightarrow \L \rightarrow \E \rightarrow \L^{-1} \rightarrow 0$$
where $\L$ is a {\bf theta characteristic}, that is, $\L ^{\otimes 2}
\cong \Omega^1 _C$.
\end{prop}

We thus have a natural set of unstable vector bundles upon which to look for connections with vanishing $p$-curvature. Indeed, it is easy to see that the proposition is sharp. 

\begin{cor} Frobenius-unstable vector bundles of rank $2$ and trivial determinant on $C$ are necessarily stable, and in one-to-one correspondence with transport-equivalence classes of connections on vector bundles $\E$ as in the above proposition, having trivial determinant and vanishing $p$-curvature. This correspondence is functorial in the sense that after arbitrary base change $C' \rightarrow C$, vector bundles $\F$ with trivial determinant and $F^* \F \cong \E'$ are in one-to-one correspondence with transport-equivalence classes of connections on $\E'$ having trivial determinant and vanishing $p$-curvature.
\end{cor}

\begin{proof} The functoriality is the more obvious statement, in light Theorem \ref{exp-cartier-vect}. For the rest, all we need check is that if $F^* \F \cong \E$ for some $\F$, we necessarily have that $\F$ is stable. But if $\M \subset \F$ is a non-negative line sub-bundle, $F^* \M \subset \F$ is non-negative with degree a multiple of $p$, which cannot occur when $F^* \F \cong \E$ by the following standard lemma.
\end{proof}

We state the lemma in more generality than immediately necessary, for later use. The argument for the $\cEnd^0(\E)$ case is taken from \cite[Lem. I.3.5, p. 105]{mo3}.

\begin{lem}\label{exp-destab-unique} Let $\E$ be a rank $2$ vector bundle of degree $0$ on a possibly nodal curve $C$, and suppose
$\L$ is a positive line bundle giving an exact sequence 
$$0 \rightarrow \L \rightarrow \E \rightarrow \E/\L
\rightarrow 0$$
Then $\L$ is unique, and is the maximal degree line bundle inside $\E$,
and $\E$ has no quotient line bundle of degree $0$. Furthermore, the same statement holds for positive sub-bundles of given rank of the traceless endomorphisms $\cEnd^0(\E)$.
\end{lem}

\begin{proof}One checks this simply by considering maps of the form
$\L\rightarrow \E \rightarrow \E/\L'$, and considering the degrees of the
line bundles in question. For the $\cEnd^0(\E)$ case, because $\cEnd^0(\E)$
is self-dual it suffices to consider the case of line sub-bundles, and to
show that the existence of a positive sub-bundle precludes the existence of
a line sub-bundle of degree $0$. But if we have $\L \subset \cEnd^0(\E)$
positive, and $\L' \subset \cEnd^0(\E)$ non-negative, first by considering
$\L' \rightarrow \cEnd^0(\E) \rightarrow \L^{-1}$ we find that the
composition must be zero, so that we have a map $\cEnd^0(\E)/\L' \rightarrow \L^{-1}$. But then considering the natural $\O_C \subset \cEnd^0(\E)/\L'$, composing with the map to $\L^{-1}$ must again give zero, so that in fact the map $\cEnd^0(\E) \rightarrow \L^{-1}$ factors through $(\cEnd^0(\E)/\L')/\O_C \cong \L'^{-1}$, from which one can conclude the desired statement.
\end{proof}

Next, we note that the $\E$ of Proposition \ref{exp-unstable} are nearly unique.

\begin{prop}\label{exp-unstable-unique} There are only $16$ choices for 
$\E$ as described in Proposition \ref{exp-unstable}, one for each choice of
$\L$.
\end{prop}

\begin{proof} Any two choices of $\L$ differ by one of the $2^{2g}=16$ line bundle of order $2$ on $C$. With $\L$ chosen, we calculate that 
$\Ext^1(\L^{-1}, \L) \cong H^0(C, \O_C) \cong k$, so the isomorphism class of $\E$ is uniquely determined.
\end{proof}

Lastly, we observe that it suffices to handle a single choice of $\E$.

\begin{cor}\label{exp-unstable-same} For any $\E, \E'$ as in Proposition 
\ref{exp-unstable-unique}, there is a canonical funtorial equivalence
between the vector bundles $\F$ of trivial determinant with $F^* \F \cong
\E$, and those with $F^* \F \cong \E'$. 
\end{cor}

\begin{proof} From Proposition \ref{exp-unstable-unique} we see that $\E$ 
and $\E'$ are related by tensoring by a $2$-torsion line bundle. The corollary is then easily verified by the bijectivity of $F^*$ on $2$-torsion line bundles.
\end{proof}

Having reduced our main question to a matter of classifying connections with vanishing $p$-curvature on a certain vector bundle, we briefly develop the formal properties of $p$-curvature, which we will not need to use until Section \ref{s-exp-det} and the following sections. The statement is:

\begin{prop}\label{exp-pcurve-formal} Given a connection $\nabla$ on a vector bundle $\E$ on a smooth $X$ over $S$, we have the following description of the $p$-curvature $\psi_{\nabla}$ of $\nabla$.
\begin{ilist}
\itm We may describe $\psi_{\nabla}$ as an element of 
$$\Gamma(X, \cEnd(\E) \otimes F^* \Omega^1_{X^{(p)}/S})^{\nabla^{\ind}},$$
where the superscript denotes the subspace of sections horizontal for
$\nabla^{\ind}$;
\itm If $\E$ and $\nabla$ have trivial determinant, we find that 
$\psi_{\nabla}$ lies in
$$\Gamma(X, \cEnd^0(\E) \otimes F^* \Omega^1_{X^{(p)}/S})^{\nabla^{\ind}},$$
where $\cEnd^0(\E)$ denotes the sheaf of traceless endomorphisms of $\E$.
\itm Assuming $\E$ has a connection, we may also consider $p$-curvature as giving maps between affine spaces
$$\psi: \Gamma(X, \Conn(\E)) \rightarrow 
\Gamma(X, \cEnd(\E) \otimes F^* \Omega^1_{X^{(p)}}),$$
$$\psi^0: \Gamma(X, \Conn^0(\E)) \rightarrow 
\Gamma(X, \cEnd^0(\E) \otimes F^* \Omega^1_{X^{(p)}}),$$
where $\Conn(\E)$ and $\Conn^0(\E)$ denote the sheaves of connections on $\E$, and of connections with trivial determinant (when $\E$ has trivialized determinant) respectively.
\itm We may take the determinant of the previous maps, and in the case that $\E$ has trivial determinant, we obtain a map
$$\det \psi^0: \Gamma(X, \Conn^0(\E)) \rightarrow 
\Gamma(X^{(p)}, (\Omega^1_{X^{(p)}})^{\otimes n}).$$
\end{ilist}
\end{prop}

\begin{proof} Assertion (i) follows directly from the linearity and $p$-linearity results of Katz \cite[5.0.5, 5.2.0]{ka1}, together with the fact that $\psi_\nabla (\theta)$ commutes with $\nabla_{\theta'}$ for any $\theta'$, by\cite[5.2.3]{ka1}. Assertion (ii) follows from explicit computation, in Corollary \ref{exp-pcurve-cors} (ii). We then obtain assertion (iii) formally: since we are working over an arbitrary scheme, we obtain the map on arbitrary $T$-valued points, and if $\E$ has a connection, the space of connections is a torsor over $\Gamma(X, \E\otimes \Omega^1_X)$, and likewise after arbitrary pull-back, and hence representable by an affine space. Finally, for assertion (iv), we just put together assertions (ii) and (iii), checking that in the trivial determinant case, the induced connection on the determinant of $\End^0(\E)$ is likewise trivial.
\end{proof}

\section{Explicit $p$-curvature Formulas}\label{s-exp-gen-pcurve}

In this section, we develop general combinatorial formulas which may be used to explicitly compute the $p$-curvature of a connection for any given $p$, and in any dimension, although it will be easiest to compute in the case of curves, where it suffices to consider a single derivation. We specify our notation for the section.

\begin{sit}\label{exp-single-open} $U$ denotes an affine open on a 
smooth $X/S$. We are given a
vector bundle $\E$ trivialized on $U$, and a derivation $\theta$ on $U$. We thus obtain a connection matrix $\bar{T}$ on $U$ associated 
to any connection $\nabla$ on $\E$, such that $\nabla _{\theta} (s) = \bar{T} s + \theta s$. Denote also by $\bar{T}_{(p)}$ the connection matrix associated to $\nabla$ and $\theta^p$.
\end{sit}

One can then easily check the following explicit formula for the $p$-curvature associated to $\nabla$ and $\theta$.

\begin{lem} We have $\psi_{\nabla} (\theta) = (\bar{T} + \theta)^p - 
\bar{T}_{(p)} - \theta ^p$ 
\end{lem}

We now describe the expansion of $(\bar{T} + \theta)^n$ using the
commutation relation $\theta \bar{T} = (\theta \bar{T}) + \bar{T} \theta$, where, in order to make formulas easier to parse, $(\theta
\bar{T})$ denotes the application of $\theta$ to the coordinates of $\bar{T}$. 

\begin{prop}Given $\i = (i_1, \dots, i_\ell) \in \N^{\ell-1}\times (\N \cup \{0\})$ 
with $\sum _{j=1} ^\ell i_j 
= n$, denote by $\hat{n}_\i$ the coefficient of $\bar{T}_{\i}:=(\theta ^{i_1-1} \bar{T})\dots 
(\theta ^{i_{\ell-1}-1} \bar{T}) \theta ^{i_\ell}$ in the full expansion of
$(\bar{T} + 
\theta)^n$. Also denote by $\i_0$ the vector $(i_1, \dots, i_{\ell-1}, 0)$. 
Then we have: 
$$\hat{n}_{\i} = \binom{n}{i_\ell} \hat{n}_{\i_0}$$
\end{prop}

\begin{proof} Although this formula may be seen directly, the proof is
expressed most clearly by induction on $n$, which we sketch. We may assume
that $i_{\ell}>0$, or the statement is trivial. By definition, we have 
$$(\bar{T}+\theta)^n= (\bar{T}+\theta)(\sum _{\ell'} \sum _{|\i'|=n-1}
\hat{n}_{\i'} \bar{T}_{\i'}),$$ 
where $\i' = (i'_1, \dots, i'_{\ell'})$, and $|\i'|:= \sum _j i'_j$.
Multiplying out and commuting the $\theta$ from left to right until we
obtain another such expression, we
find two cases: $i_1 =1$ and $i_1>1$; we handle the case $i_1>1$,
the other being essentially the same. In this case, we
obtain the inductive formula $\hat{n}_{\i} = \sum_j \hat{n}_{\i - 1_j}$,
where $1_j$ denotes the vector which is $1$ in the $j$th position and $0$
elsewhere, and where $j$ is allowed to range only over values where $i_j >
1$. We then have also that $\hat{n}_{\i_0} = \sum_{j<\ell} \hat{n}_{\i_0-1_j}$,
so that if we induct on $n$, we have 
$\hat{n}_i = \sum _j \hat{n}_{i-1_j} = 
\sum _{j <\ell} \binom{n-1}{i_\ell} \hat{n}_{(\i-1_j)_0} + \binom{n-1}{i_{\ell-1}}
\hat{n}_{(\i - 1_\ell)_0} = (\binom{n-1}{i_\ell}+\binom{n-1}{i_\ell-1})
\hat{n}_{\i_0},$ where the last equality makes use of the observation that
$\i_0 = (\i-1_\ell)_0$. Then the identity $\binom{n-1}{r}+ \binom{n-1}{r-1} =
\binom{n}{r}$ completes the proof.
\end{proof}

It follows that if $n=p$, $\hat{n}_{\i}$ is nonzero mod $p$ only if $i_\ell=0$ or
$i_\ell =p$, and in the latter case, we have $\ell=1$, $i_1 =p$, and 
$\hat{n}_{\i} =1$, 
which precisely cancels the $\theta ^p$ subtracted off in the formula
for $\psi_{\nabla}(\theta)$. We immediately
see that $\psi_{\nabla}(\theta)$ is in fact given entirely by linear
terms. In particular, this explicitly recovers the statement we already knew to be true that
$p$-curvature takes values in the space of $\O_C$-linear endomorphisms of 
$\E$. We may now restrict our attention to the linear terms in the 
expansion, and will shift our notation accordingly:

\begin{prop}\label{exp-gen-pcurve} Given $\i = (i_1, \dots, i_\ell) \in \N^\ell$ with $\sum _{j=1}
^\ell i_j 
= n$, denote by $n_\i$ the coefficient of $\bar{T}_{\i}=(\theta ^{i_1-1} \bar{T})\dots 
(\theta ^{i_\ell-1} \bar{T})$ in the full expansion of $(\bar{T} + \theta)^n$. Also
denote by $\hat{\i}$ the truncated vector $(i_1, \dots, i_{\ell-1})$. Then we 
have: 
$$n_{\i} = \binom{n-1}{i_\ell-1} n_{\hat{\i}}$$
We thus get 
$$n_{\i} = \prod _{j=1} ^\ell \binom{n-1 - \sum _{m=j+1} ^\ell i_m}{i_j -1} 
= \frac{(n-1)!}{(\prod_{j=1}^\ell (i_j-1)!)(\prod _{j=1} ^{\ell-1} (\sum _{m=1}
^j i_m))}$$
\end{prop}

\begin{proof} This follows from the same induction argument as the previous proposition.
\end{proof}

We note that this implies that every such term in the expansion
of $(\bar{T}+\theta)^n$ is nonzero mod $n$ when $n=p$, since the numerator in the
resulting formula is simply $(n-1)!$. Thus, the $p$-curvature formula is 
always maximally complex, having an exponential number of terms. However, when some of the terms commute, the
formulas tend to simplify considerably.

\begin{prop}\label{exp-pcurve-commute} Given $\ell>0$ and a subset 
$\Lambda \subset \{1, \dots, \ell\}$, denote by $S^{\Lambda}_{\ell}$ 
the subset of the permutation group $S_{\ell}$ which preserves the order of
the elements of $\Lambda$; that is, $S^{\Lambda}_{\ell}:=\{\sigma \in
S_{\ell}:\forall j<j' \in \Lambda, \sigma(j)<\sigma(j')\}$. Given also
$\i = (i_1, \dots, i_\ell) \in \N^\ell$ with $\sum _{j=1} ^\ell i_j 
= n$, we denote by $n^{\Lambda}_\i$ the 
sum over all $\sigma \in S^\Lambda_{\ell}$ of $n_{\sigma(\i)}$, where
$\sigma(\i)$ denotes the vector $(i_{\sigma^{-1}(1)}, \dots,
i_{\sigma^{-1}(\ell)})$ obtained from $\i$ by permuting the coordinates
under $\sigma$. Then we have: 
$$n^{\Sigma}_{\i} = \frac{n!}{ \prod _{j=1} ^\ell (i_j-1)!
\prod _{j=1}^{\ell}(i_j+\sum_{m<j}^{m,j \in \Lambda} i_m)}.$$
Note that the last
sum in the denominator is non-empty only for $j \in \Lambda$.
\end{prop}

\begin{proof} First note that if we want the entries of $\i$ with indices in
$\Lambda$ to have the same order in $\sigma(\i)$, we must apply
$\sigma^{-1}$ rather than $\sigma$ to the indices, as in our definition. 
Applying our previous formula, we really just want to show that 
$$\sum _{\sigma \in S^{\Lambda}_\ell} \prod _{j=1} ^{\ell-1} 
\frac{1}{\sum _{m=1} ^j i_ {\sigma^{-1}(m)}}
= \frac{n}{\prod _{j=1} ^\ell (i_j + \sum_{m<j}^{m,j \in \Lambda} i_m)} = 
\frac{\sum_{j=1} ^\ell i_j}{\prod _{j=1} ^\ell (i_j + \sum_{m<j}^{m,j \in
\Lambda} i_m)}.$$
Dividing through by $\sum_{j=1}^\ell i_j$ reduces the identity to
\begin{equation}\label{exp-perm-comb}
\sum _{\sigma \in S^{\Lambda}_\ell} \prod _{j=1} ^{\ell} 
\frac{1}{\sum _{m=1} ^j i_ {\sigma^{-1}(m)}} =
\frac{1}{\prod _{j=1} ^\ell (i_j + \sum_{m<j}^{m,j \in \Lambda} i_m)}.
\end{equation}
We show this by induction on $\ell$ (noting that it is rather trivial in the
case $\ell=1$, whether or not $\Lambda$ is empty), breaking up the first sum 
over $S^{\Lambda}_\ell$ into $\ell-|\Lambda|+1$ pieces, depending on which 
$i_r$ ends up in the final place. There are two cases to consider: 
$r \not\in \Lambda$,
or $r = \Lambda_{\max}$. In either case, the relevant part of the sum on 
the left hand side becomes $\sum _{\sigma \in S^{\Lambda,r}_{\ell}} 
\prod _{j =1} ^{\ell} \frac{1}{\sum _{m=1}^j i_{\sigma^{-1}(m)}}$, where
$S^{\Lambda,r}_{\ell}$ denotes the subset of $S^{\Lambda}_{\ell}$ sending
$r$ to $\ell$. Now, the point is that for our sums, this will be equivalent
to an order-preserving subset of the symmetric group acting on a set of 
$\ell-1$ elements, allowing us to apply induction. In 
the case that $r \not \in \Lambda$, $\Lambda$ is in essence unaffected, and
we find that 
$$\sum _{\sigma \in S^{\Lambda,r}_{\ell}} \prod _{j =1} ^{\ell} 
\frac{1}{\sum _{m=1}^j i_{\sigma^{-1}(m)}}= \frac{1}{n} 
\sum _{\sigma \in S^{\Lambda,r}_{\ell}} \prod _{j = 1}^{\ell-1} 
\frac{1}{\sum _{m=1}^{j} i_{\sigma^{-1}(m)}},$$ 
and one checks that this sum is of the same form as Equation
\ref{exp-perm-comb}, with $i_r$ ommitted, so by induction we find that this
sum is equal to
$$ \frac{1}{n}\frac{1}{\prod_{j \neq r} (i_j + \sum_{m<j}^{m,j \in \Lambda}
i_m)} = \frac{i_r+\sum_{m<r}^{m,r \in \Lambda}i_m}{n\prod_{j=1}^\ell 
(i_j + \sum_{m<j}^{m,j \in \Lambda} i_m)}
= \frac{i_r}{n \prod_{j=1}^\ell (i_j + \sum_{m<j}^{m,j \in \Lambda} 
i_m)},$$ 
since $r \not \in \Lambda$.
In the case that $r=\Lambda_{\max}$, we effectively reduce the size
of $\Lambda$ by one, but because $r$ is maximal in $\Lambda$, for $j \neq r$
the term $\sum _{m <j}^{m,j \in \Lambda}i_m $ is unaffected by omitting $r$ from
$\Lambda$. We thus find, arguing as before,
$$\sum _{\sigma \in S^{\Lambda,r}_{\ell}} \prod _{j =1} ^{\ell} 
\frac{1}{\sum _{m=1}^j i_{\sigma^{-1}(m)}} = 
\frac{i_r+\sum_{m<r}^{m,r \in \Lambda}i_m}{n \prod_{j=1}^\ell 
(i_j + \sum_{m<j}^{m,j \in \Lambda} i_m)}
= \frac{\sum_{j \in \Lambda} i_j}{n \prod_{j=1}^\ell 
(i_j + \sum_{m<j}^{m,j \in \Lambda} i_m)}.$$ 
Adding these up as $r$ ranges over $\Lambda_{\max}$ and all values not in
$\Lambda$, and using $n = \sum _j i_j$, we get the desired identity.
\end{proof}

We give some specific applications of this formula.

\begin{cor}\label{exp-pcurve-cors} Let $\E$ be a vector bundle of rank $r$ on a smooth variety $X$ over a field $k$, with $\nabla$ an integrable connection on $\E$ and $\theta$ a derivation on an open set $U$ which also trivializes $\E$. We have:
\begin{ilist}
\itm If $r=1$, $p$-curvature is given by:
$$\psi_{\nabla}(\theta) = 
\bar{T}^p + (\theta ^{p-1} \bar{T}) - \bar{T}_{(p)}.$$
\itm Suppose $\E$ has trivialized determinant, and $\nabla$ has trivial determinant. Then the $p$-curvature of $\nabla$ has image in the traceless endomorphisms of $\E$.
\itm Suppose $\nabla'$ is a connection on $U$ with $\nabla'-\nabla = \omega \Id$ a scalar
endomorphism. Then we have
$$\psi_{\nabla'}(\theta)-\psi_{\nabla}(\theta) = ((\hat{\theta}(\omega))^p +
\theta^{p-1}(\hat{\theta}(\omega))-\hat{\theta^p}(\omega))\Id,$$ 
where $\hat{\theta}$ denotes the unique linear map $\Omega^1_C\rightarrow
\O_C$ such that $\theta = \hat{\theta} \circ d$.
\end{ilist}
\end{cor}

\begin{proof} From the previous proposition, we see that when $n=p$ and $\Lambda$ is empty, so that all the involved matrices commute, we have 
$$n^{\varnothing}_{\i} = \frac{p!}{\prod_{j=1}^{\ell}i_j!},$$ 
but the actual coefficient will be $n^{\varnothing}_{\i}/P_{\i}$, where
$P_{\i}$ is the number of permutations fixing the vector $\i$, since summing
up over all permutations will count each term $P_{\i}$ times.
We see that this expression can be non-zero mod $p$
only if either $P_{\i}$ is a multiple of $p$, or some $i_j$ is. Since $P_{\i}$
is the order of a subgroup of $S_\ell$, it can be a multiple of $p$ if and only
if $\ell=p$ and each $i_j=1$. On the other hand, an $i_j$ can be a multiple
of $p$ if and only if $\ell=1$ and $i_1=p$; these two terms simply reiterate
that the coefficients of $\bar{T}^p$ and $(\theta ^{p-1} \bar{T})$ are both 
$1$, and we see that every other coefficient vanishes mod $p$. We immediately conclude (i), and for (ii) we see similarly that we have
$$\Tr \psi_{\nabla}(\theta) = 
\Tr \bar{T}^p + \Tr (\theta ^{p-1} \bar{T}) - \Tr f_{\theta ^p} \bar{T}.$$
The second and third terms visibly have vanishing trace because $\bar{T}$
does, while it is easy to see (for instance, by passing to the algebraic
closure of $k$ and taking the Jordan normal form) that $\Tr (\bar{T}^p) =
(\Tr \bar{T})^p =0$.

For (iii), we can compare the $p$-curvatures of $\nabla$ and $\nabla'$ term 
by term; we have $\bar{T}+\hat{\theta}(\omega) \Id$ as the matrix for 
$\nabla'$, and we see that if we expand each term of 
$\psi_{\nabla'}(\theta)$, we get $\psi_{\nabla}(\theta)$ from expanding out 
only terms involving $\bar{T}$ and $\theta$, and $((\hat{\theta}(\omega))^p 
+ \theta^{p-1}(\hat{\theta}(\omega))-\hat{\theta^p}(\omega))\Id$ 
from expanding out terms involving only $\hat{\theta}(\omega)$ and 
$\theta$, since these last all commute with one another. We thus want to 
show that all of the coefficients of the cross terms are always zero mod 
$p$. If we consider a particular term 
$(\theta_0 ^{i_1-1} (\bar{T}+\hat{\theta}(\omega)\Id))\dots 
(\theta_0 ^{i_\ell-1} (\bar{T}+\hat{\theta}(\omega)\Id))$
corresponding to a vector $\i$, a cross term will arise by choosing a 
subset $\Lambda \subset \{1, \dots, \ell\}$ from which the $\bar{T}$ terms 
will be chosen, with the $\hat{\theta}(\omega)\Id$ term being chosen for 
all indices outside $\Lambda$. To compute the relevant coefficient we can 
essentially sum over all permutations in the $S^{\Lambda}_{\ell}$ of
Proposition \ref{exp-pcurve-commute}. The only caveat is that if 
$\sigma \in S^{\Lambda}_{\ell}$ fixes $\Lambda$ and leaves the vector $\i$ unchanged, then it will give the same term in the expansion as the 
identity. Such $\sigma$ form a subgroup of $S_{\ell}$, and if we denote the 
order of this subgroup by $P^{\Lambda}_{\i}$, we find that the coefficient 
we want to compute is given by, still in the notation of Proposition 
\ref{exp-pcurve-commute}, the expression 
$n^{\Lambda}_{\i}/P^{\Lambda}_{\i}$. Now, the only way to cancel the $p$ in
the numerator of $n^{\Lambda}_{\i}$ would be for either $P^{\Lambda}_{\i}$ 
or the denominator of $n^{\Lambda}_{\i}$ to also be divisible by $p$. The 
denominator of $n^{\Lambda}_{\i}$ cannot be divisible by $p$, since the 
$i_j$ add up to $p$, and the only way that $p$ could appear in the 
denominator would therefore be when $\Lambda$ is all of 
$\{1, \dots, \ell\}$, which corresponds to the terms which only involve 
$\bar{T}$, or when $\ell=1$, which gives the 
$\theta^{p-1}(\hat{\theta}(\omega))$ term. Similarly, $P^{\Lambda}_{\i}$ is 
the order of a subgroup of $S_{\ell}$ which fixes $\Lambda$, so can be a 
multiple of $p$ only if $\ell=p$ and $|\Lambda|=0$, which corresponds to 
the term $((\hat{\theta}(\omega))^p$. This yields the desired result.
\end{proof}

We do not use the last statement of the corollary in this paper, but it could be used to show, for instance, that when $r$ is prime to $p$, and $\E$ has rank $r$, then if any representative on $\E$ of a projective connection $\nabla$ on $\P(\E)$ has vanishing $p$-curvature, then the unique representative on $\E$ of $\nabla$ with vanishing trace must likewise have vanishing $p$-curvature. We also remark that results such as statement (ii) above may generally be obtained more abstractly via general functoriality statements on $p$-curvature, but such a point of view requires familiarity with Grothendieck's abstract theory of connections; see \cite{os6}.

We conclude with some observations in the case of curves. 

\begin{lem} In the case that $X$ is a curve, the $p$-curvature of a connection $\nabla$ is identically $0$ if and only if $\psi_{\nabla} (\theta)=0$ for any non-zero derivation 
$\theta$. In addition, $\bar{T}_{(p)}=f_{\theta^p} \bar{T}$ for some
function $f_{\theta^p}$, satisfying $f_{\theta^p}\theta = \theta^p$.
\end{lem}

\begin{proof} These statements follow trivially from the fact that the sheaf of derivations is invertible, and the $p$-linearity of the $p$-curvature map with respect to derivations.
\end{proof}

Finally, we record in this situation the general $p$-curvature formulas in 
characteristics $3$, $5$, and $7$, for later use.

Characteristic 3:
\begin{equation}\label{exp-pcurve-3}
\psi_\nabla (\theta)=
\bar{T}^3+(\theta \bar{T})\bar{T} + 2 \bar{T} (\theta  \bar{T}) +(\theta^2 \bar{T})-f_{\theta ^3} \bar{T}
\end{equation}

Characteristic 5:
\begin{multline}\label{exp-pcurve-5}
\psi_\nabla (\theta)= 
\bar{T} ^5 + 4 \bar{T} ^3(\theta^1 \bar{T}) + 3 \bar{T} ^2(\theta^1 \bar{T})\bar{T} + \bar{T} ^2(\theta^2 \bar{T}) + 
2 \bar{T}(\theta^1 \bar{T})\bar{T} ^2 \\
+ 3 \bar{T}(\theta^1 \bar{T}) ^2 + 3 \bar{T}(\theta^2 \bar{T})\bar{T} + 4 \bar{T}(\theta^3 \bar{T}) + 
(\theta^1 \bar{T})\bar{T} ^3 \\
+ 4 (\theta^1 \bar{T})\bar{T}(\theta^1 \bar{T}) 
+ 3 (\theta^1 \bar{T}) ^2\bar{T} + (\theta^1 \bar{T})(\theta^2 \bar{T}) +
(\theta^2 \bar{T})\bar{T} ^2 \\
+ 4 (\theta^2 \bar{T})(\theta^1 \bar{T}) + (\theta^3 \bar{T})\bar{T} 
+ (\theta^4 \bar{T}) - f_{\theta^5} \bar{T}
\end{multline}

Characteristic 7:
\begin{multline}\label{exp-pcurve-7}
\psi_\nabla (\theta)=
\bar{T} ^7 + 6 \bar{T} ^5(\theta^1 \bar{T}) + 5 \bar{T} ^4(\theta^1 \bar{T})\bar{T} + \bar{T} ^4(\theta^2 \bar{T}) +
4 \bar{T} ^3(\theta^1 \bar{T})\bar{T} ^2 \\
+ 3 \bar{T} ^3(\theta^1 \bar{T}) ^2 + 3 \bar{T} ^3(\theta^2 \bar{T})\bar{T} +
6 \bar{T} ^3(\theta^3 \bar{T}) + 3 \bar{T} ^2(\theta^1 \bar{T})\bar{T} ^3 \\
+ 4 \bar{T} ^2(\theta^1 \bar{T})\bar{T}(\theta^1 \bar{T}) + \bar{T} ^2(\theta^1 \bar{T}) ^2\bar{T} + 
3 \bar{T} ^2(\theta^1 \bar{T})(\theta^2 \bar{T}) + 6 \bar{T} ^2(\theta^2 \bar{T})\bar{T} ^2 \\
+ \bar{T} ^2(\theta^2 \bar{T})(\theta^1 \bar{T}) + 3 \bar{T} ^2(\theta^3 \bar{T})\bar{T} +
\bar{T} ^2(\theta^4 \bar{T}) + 2 \bar{T}(\theta^1 \bar{T})\bar{T} ^4 \\
+ 5 \bar{T}(\theta^1 \bar{T})\bar{T} ^2(\theta^1 \bar{T}) + 3 \bar{T}(\theta^1 \bar{T})\bar{T}(\theta^1 \bar{T})\bar{T} + 
2 \bar{T}(\theta^1 \bar{T})\bar{T}(\theta^2 \bar{T}) \\
+ \bar{T}(\theta^1 \bar{T}) ^2\bar{T} ^2 + 6 \bar{T}(\theta^1 \bar{T}) ^3 + 
6 \bar{T}(\theta^1 \bar{T})(\theta^2 \bar{T})\bar{T} + 5 \bar{T}(\theta^1 \bar{T})(\theta^3 \bar{T}) \\
+ 3 \bar{T}(\theta^2 \bar{T})\bar{T} ^3 + 4 \bar{T}(\theta^2 \bar{T})\bar{T}(\theta^1 \bar{T}) + 
\bar{T}(\theta^2 \bar{T})(\theta^1 \bar{T})\bar{T} + 3 \bar{T}(\theta^2 \bar{T}) ^2 \\
+ 4 \bar{T}(\theta^3 \bar{T})\bar{T} ^2 + 3 \bar{T}(\theta^3 \bar{T})(\theta^1 \bar{T}) + 5 \bar{T}(\theta^4 \bar{T})\bar{T} + 
6 \bar{T}(\theta^5 \bar{T}) + (\theta^1 \bar{T})\bar{T} ^5 \\
+ 6 (\theta^1 \bar{T})\bar{T} ^3(\theta^1 \bar{T}) + 5 (\theta^1 \bar{T})\bar{T} ^2(\theta^1 \bar{T})\bar{T} + 
(\theta^1 \bar{T})\bar{T} ^2(\theta^2 \bar{T}) \\
+ 4 (\theta^1 \bar{T})\bar{T}(\theta^1 \bar{T})\bar{T} ^2 + 3 (\theta^1 \bar{T})\bar{T}(\theta^1 \bar{T}) ^2  + 
3 (\theta^1 \bar{T})\bar{T}(\theta^2 \bar{T})\bar{T} \\
+ 6 (\theta^1 \bar{T})\bar{T}(\theta^3 \bar{T}) + 3 (\theta^1 \bar{T}) ^2\bar{T} ^3 + 
4 (\theta^1 \bar{T}) ^2\bar{T}(\theta^1 \bar{T}) + (\theta^1 \bar{T})
^3\bar{T} \\
+ 3 (\theta^1 \bar{T}) ^2(\theta^2 \bar{T}) + 6 (\theta^1 \bar{T})(\theta^2 \bar{T})\bar{T} ^2 + 
(\theta^1 \bar{T})(\theta^2 \bar{T})(\theta^1 \bar{T}) \\
+ 3 (\theta^1 \bar{T})(\theta^3 \bar{T})\bar{T} + (\theta^1 \bar{T})(\theta^4 \bar{T}) 
+ (\theta^2 \bar{T})\bar{T} ^4 + 6 (\theta^2 \bar{T})\bar{T} ^2(\theta^1 \bar{T}) \\
+ 5 (\theta^2 \bar{T})\bar{T}(\theta^1 \bar{T})\bar{T} + (\theta^2 \bar{T})\bar{T}(\theta^2 \bar{T}) + 
4 (\theta^2 \bar{T})(\theta^1 \bar{T})\bar{T} ^2 \\
+ 3 (\theta^2 \bar{T})(\theta^1 \bar{T}) ^2 + 3 (\theta^2 \bar{T}) ^2\bar{T} + 
6 (\theta^2 \bar{T})(\theta^3 \bar{T}) + (\theta^3 \bar{T})\bar{T} ^3 \\ 
+ 6 (\theta^3 \bar{T})\bar{T}(\theta^1 \bar{T}) + 5 (\theta^3 \bar{T})(\theta^1 \bar{T})\bar{T} + 
(\theta^3 \bar{T})(\theta^2 \bar{T}) + (\theta^4 \bar{T})\bar{T} ^2 \\
+ 6 (\theta^4 \bar{T})(\theta^1 \bar{T}) + (\theta^5 \bar{T})\bar{T} + (\theta^6 \bar{T}) - 
f_{\theta^7} \bar{T}
\end{multline}

\section{On $f_{\theta^p}$ and $p$-rank in Genus $2$}\label{s-exp-ftheta}

In this section, we give an explicit formula for $f_{\theta^p}$ on a 
genus $2$ curve $C$, and note that we can use these ideas to derive explicit
formulas for the $p$-rank of the Jacobian of $C$. Throughout, we work under
the hypotheses and notation of Situations \ref{exp-genus} and
\ref{exp-single-open}, with $X=C$.

We first note that (irrespective of the genus of $C$), although $f_{\theta ^p}$ will be $0$ only if $\theta(f)=1$ for some $f$ on $U$, we will always have:

\begin{lem}\label{exp-thetaftheta} $\theta f_{\theta ^p}=0$.
\end{lem}

\begin{proof} Given any $f$, $\theta ^p f = f_{\theta^p}\theta(f)$, so
$\theta ^{p+1} f = \theta (f_{\theta ^p} \theta(f)) =
\theta (f_{\theta ^p}) \theta(f) + f_{\theta ^p} \theta^2(f)
= \theta (f_{\theta ^p}) \theta(f) + \theta ^{p+1}(f)$.
Since this is true for all $f$, we must have $\theta(f_{\theta^p})=0$, as
desired.
\end{proof}

We now specify some normalizations and notational conventions special to
genus $2$ which we will follow through the end of our explicit calculations in Section \ref{s-exp-det}.

\begin{sit}\label{exp-presented} $C$ is a smooth, proper genus $2$
curve over an algebraically closed field $k$. It is presented explicitly on
an affine open set $U_2$ by 
$$y^2 = g(x) = x^5 +a_1 x^4 + a_2 x^3 + a_3 x^2 + a_4 x + a_5,$$ with the
complement of $U_2$ being a single, smooth, Weierstrass point $w$ at 
infinity. We also have the form $\omega_2 = y^{-1} dx$ trivializing
$\Omega^1_C$ on $U_2$, and the derivation $\theta$ on $U_2$ given by
$\theta f = y \frac{df}{dx}$. Equivalently, $\hat{\theta}(\omega_2)=1$,
where $\hat{\theta}$ denotes the map $\Omega^1_C \rightarrow \O_C$ such that
$\theta= \hat{\theta} \circ d$.
\end{sit}

For this section only, we set $U=U_2$ and $\omega=\omega_2$. 
We set $g_k(x) = \theta ^{k-1} x$;
we see by induction that this is a polynomial in $x$ for $k$ odd. 
Noting that $\theta(p(x))=yp'(x)$ for $p(x)$ any polynomial in $x$, and 
$\theta(y)=\frac{1}{2}g'(x)$, we have
that for $k$ odd, $g_k(x) = \theta^2 (g_{k-2}(x))= \theta (y
g'_{k-2}(x))$, and we get the recursive formula:
\begin{equation}\label{exp-gk} g_k(x)=g''_{k-2}(x)g(x) + 
\frac{1}{2}g'_{k-2}(x)g'(x) \end{equation}
for $k$ odd.

But $f_{\theta^p} = \hat{\theta}^p(y^{-1} dx)$ by definition, which is
just $y^{-1} \theta ^p (x)$, so we 
also find 
\begin{equation}\label{exp-ftheta} f_{\theta^p} = 
y^{-1} \theta g_p (x) = g'_p(x) \end{equation}
In particular, $f_{\theta^p}$ is a polynomial in $x$, and can therefore
only have nonzero terms mod $p$ in degrees which are multiples of $p$.
However, we see by induction that the degree of 
$g_p(x)$ is always less than $2p$, so the only nonzero terms of
$f_{\theta^p}$ are the constant term and the $p$th power term (from
which it follows that the only nonzero terms of $g_p(x)$ are the constant,
linear, $p$th power, and $(p+1)$st power terms).

For later use, we note the formulas for characteristics $3$, $5$, and $7$
obtained by combining equations \ref{exp-gk} and \ref{exp-ftheta}:

Characteristic 3:
\begin{equation}\label{exp-ftheta-3}f_{\theta^3} =
x^3+a_3 \end{equation}

Characteristic 5:
\begin{equation}\label{exp-ftheta-5}f_{\theta^5} = 
2 a_1 x^5 + a_3^2 + 2 a_2 a_4 + 2 a_1 a_5 \end{equation}

Characteristic 7:
\begin{equation}\label{exp-ftheta-7}f_{\theta^7} = 
(3 a_1 ^2 + 3 a_2) x ^7 + a_3 ^3 + 6 a_2 a_3 a_4 + 3 a_1 a_4 ^2 + 
3 a2 ^2 a_5 + 6 a_1 a_3 a_5 + 6 a_4 a_5 \end{equation}

As a final note, we can use this to derive explicit formulas for the
$p$-rank of the Jacobian of $C$ in terms of the coefficients of $g(x)$.

\begin{prop}\label{exp-prank} If we denote by $h_1, h_2, h_3, h_4$ the 
polynomials in the coefficients of $g(x)$ giving the constant, linear, 
$p$th power, and $(p+1)$st power terms of $g_p(x)$, then the $p$-rank of 
the Jacobian of $C$ is:
\begin{itemize}
\item[$2$ if:] $h_1 h_4- h_2 h_3 \neq 0$;
\item[$1$ if:] $h_1 h_4- h_2 h_3 =0$ but either $h_3^p - h_2 h_4 ^{p-1}
\neq 0$ or $h_1^p h_4 - h_2 ^{p+1} \neq 0$;
\item[$0$ if:] $h_1 h_4- h_2 h_3 = h_3^p - h_2 h_4 ^{p-1} = h_1^p h_4 -
h_2^{p+1} = 0$.
\end{itemize}
\end{prop}

\begin{proof} The $p$-torsion of the Jacobian is simply the number of
(transport equivalence classes of, but endomorphisms of a line bundle
are only scalars, and hold connections fixed) connections with 
$p$-curvature $0$ on the trivial bundle. We note that the space of
connections on $\O_C$ can be written explicitly as $f \mapsto df + f (c_1 
+c_2 x) \omega$, meaning the connection matrix on $U$ with respect to $\theta$ is given simply by 
the function $\bar{T}= c_1+c_2 x$. Using the $p$-curvature formula for 
rank $1$ given by Corollary \ref{exp-pcurve-cors} (i), we find
\begin{gather} \psi_{\nabla}(\theta_0) = (c_1 + c_2 x)^p + 
\theta_0^{p-1}(c_1+c_2 x) - f_{\theta_0 ^p} (c_1 + c_2 x) \\
= c_1 ^p + c_2 ^p x^p + c_2 g_p(x) - g'_p(x) (c_1 +c_2 x) \\
= (c_1^p + c_2 h_1 - c_1 h_2) + (c_2 ^p + c_2 h_3 - c_1 h_4) x^p.
\end{gather}
Setting the $p$-curvature to zero, we obtain:
$$0 = (c_1^p + c_2 h_1 - c_1 h_2) + (c_2 ^p + c_2 h_3 - c_1 h_4) x^p.$$

We first consider this equation in the case that $h_4 \neq 0$.
In this case, we find that we can write $c_1= \frac{c_2^p+c_2 h_3}{h_4}$, and substituting in, we find we get $p^2$ solutions if
$h_1 h_4^p - h_2 h_3 h_4^{p-1} \neq 0$, and otherwise, $p$ solutions if 
$h_3 ^p - h_2 h_4^{p-1} \neq 0$, and finally $1$ solution if both vanish.
On the other hand, in the case that $h_4=0$, we see that $c_2$ becomes independent of $c_1$, we get $p^2$ solutions if and only if both $h_2$ and 
$h_3$ are nonzero; $p$ solutions if either but not both are nonzero, and 
$1$ solution if they are both $0$. One can then check that both these casese are expressed by the asserted polynomial conditions in the $h_i$.
\end{proof}

For $p=3$, we have 
$$g_p(x) = 1 x^4 - a_1 x^3 + a_3 x - a_4,$$ 
so $h_4$ is always nonzero, and we find that the $p$-rank of $C$ is $2$ when 
$a_4 - a_1 a_3 \neq 0$, is $1$ when $a_4 - a_1 a_3 = 0$ but 
$a_1 ^3 - a_3 \neq 0$, and is $0$ when $a_4 - a_1 a_3 = a_1 ^3 - a_3 = 0$.

For $p=5$, we have 
$$g_p(x) = 2 a_1 x_6 + (4 a_1 ^2 + 3 a_2) x^5 + 
(a_3 ^2 + 2 a_2 a_4 + 2 a_1 a_5) x +(3 a_3 a_4 + 3 a_2 a_5),$$
so the $p$-rank of $C$ is $2$ when 
$$a_1 (a_3 a_4 + a_2 a_5) 
- (4 a_1 ^2 + 3 a_2)(a_3^2 + 2 a_2 a_4 + 2 a_1 a_5) \neq 0.$$ 
The $p$-rank is $1$ when 
$$a_1 (a_3 a_4 + a_2 a_5) 
- (4 a_1 ^2 + 3 a_2)(a_3^2 + 2 a_2 a_4 + 2 a_1 a_5) = 0$$
but either 
$$4 a_1 ^{10} + 3 a_2 ^5 - (a_3 ^2 + 2 a_2 a_4 + 2 a_1 a_5) a_1^4 \neq 0$$
or 
$$(3 a_3 ^5 a_4 ^5 + 3 a_2^5 a_5^5) 2 a_1 - (a_3^2 + 2 a_2 a_4 + 2 a_1
a_5)^6 \neq 0.$$
Lastly, the $p$-rank is $0$ when 
\begin{align*}0= a_1 (a_3 a_4 + a_2 a_5) - (4 a_1 ^2 + 3 a_2)(a_3^2 
+ 2 a_2 a_4 + 2 a_1 a_5) \notag \\
=  4 a_1 ^{10} + 3 a_2 ^5 - (a_3 ^2 + 2 a_2 a_4 + 2 a_1 a_5) a_1^4 \\
= (3 a_3 ^5 a_4 ^5 + 3 a_2^5 a_5^5) 2 a_1 - (a_3^2 + 2 a_2 a_4 + 2 a_1
a_5)^6.\notag
\end{align*}

While explicit computations of the $p$-rank of the Jacobian of a curve are
not hard in general, it is perhaps worth mentioning that this method, aside
from providing a complete and explicit solution for genus $2$ curves, does
so in a sufficiently elementary way that it can be presented as a 
calculation of the $p$-torsion of $\Pic(C)$ without knowing any properties
of the Jacobian, or even that it exists.

\section{The Space of Connections}\label{s-exp-conns}

In this section we carry out the first portion of the necessary computations
for the explicit portion of Theorem \ref{exp-main}, by calculating the space of 
transport-equivalence classes of connections on a particular vector bundle
$\E$. We suppose: 

\begin{sit}\label{exp-specific-e} With the notation and hypotheses of Situation
\ref{exp-presented}, we further declare that $\E$ is the bundle determined by 
Propositions \ref{exp-unstable} and \ref{exp-unstable-unique} for the choice 
$\L = \O_C([w])$.
\end{sit}

In this situation, if $U_1, U_2$ 
are a trivializing cover for $\L$, with
transition function $\vp _{12}$, then $\L^{-1}, \L^{\otimes 2} = \Omega^1_C,$
and $\E$ are all trivialized by this cover as well, and $\E$ can be
represented with a transition matrix of the form 
$$E=\begin{bmatrix}
\vp _{12} & \vp _\E \\
0 & \vp _{12} ^{-1}
\end{bmatrix}$$
for some $\vp _\E$ regular on $U_1 \cap U_2$.

We see immediately that we can choose $\vp_{12}$ and $U_1$ so that $\vp_{12}$ is
regular on $U_1$ with a simple zero at $w$, and non-vanishing elsewhere:
we simply set $\vp_{12}$ to be any function with a simple zero at $w$, 
and take $U_1$ to be the complement of any other zeroes and poles. For
compatibility of trivializations of $\L$ and $\Omega_C^1$, we must then
set $\omega_1 = \vp_{12}^{-2} \omega _2$. Beyond these properties, our
specific choice of $\vp_{12}$ will be completely irrelevant, but we note
that it is possible to choose $\vp_{12}$ to vary algebraically (in fact, to
be in some sense invariant) as our $a_i$ and the corresponding curves vary: we can simply
set $\vp_{12} = \frac{x^2}{y}$.

\begin{prop}\label{exp-extension} The unique non-trivial isomorphism class for $\E$ may be
realized by setting $\vp _\E = \vp_{12} ^{-2}$. 
\end{prop}

\begin{proof}We claim that there cannot be a splitting map from $\E$ back to $\L$. Indeed, one checks explicitly that such a splitting would require the existence of a rational function on $C$ having a pole of order exactly $3$ at $w$, and regular elsewhere, which is not possible.
\end{proof}

We now note that since $\vp_{12}$ has a simple zero at $w$, and $\omega_1$ is invertible at $w$, if we further restrict $U_1$ we can guarantee that $\frac{d\vp_{12}}{\omega_1}$ is likewise everywhere invertible
on $U_1$. Having done so, $\vp_\E = \vp_{12}^{-2}$, so $d\vp_\E=
-2\vp_{12}^{-3} d\vp_{12}$, and $\frac{d\vp_\E}{\omega_1}$ is regular
and nonvanishing on $U_1$ except for a pole of order $3$ at $w$.

Now, we can trivialize $\E \otimes \Omega^1_C$ on the $U_i$ with 
transition matrix $\vp _{12} ^2 E$. We can then represent a 
connection $\nabla: \E \rightarrow \E \otimes \Omega^1_C$ by 
$2 \times 2$ connection matrices $\bar{T}_1$ and $\bar{T}_2$ of functions
regular on $U_1$ and $U_2$ respectively.
These act by sending $s_i \mapsto \bar{T}_i s_i + \frac{ds_i}{\omega_i}$ on
$U_i$, where the $s_i$ are given as vectors under the trivialization,
so one checks that $\bar{T}_1$ and $\bar{T}_2$ must be related by:
$$\bar{T}_1 = \vp _{12}^2 E \bar{T}_2 E^{-1}+E \frac{dE^{-1}}{\omega _1}$$

We now explicitly compute $\bar{T}_2$ in terms of $\bar{T}_1$ in preparation for
computing the space of connections. If $\bar{T}_2=
\begin{bmatrix}
f_{11} & f_{12} \\
f_{21} & f_{22}
\end{bmatrix}$, then:

\begin{equation}\label{exp-conneq}
\bar{T}_1  = 
\begin{bmatrix}
\vp _{12}^2f_{11}+f_T &
\vp_{12}^4 f_{12} + \vp_{12}^3 \vp_\E (f_{22} - f_{11})
- \vp _{12}\vp_\E f_T -
\vp_{12}\frac{d\vp_\E}{\omega_1} \\
f_{21} & 
\vp _{12}^2f_{22}-f_T
\end{bmatrix}
\end{equation}
where $f_T=\vp_{12} \vp_\E f_{21} - \vp_{12}^{-1} \frac{d \vp_{12}}{\omega_1}$

Note that this implies $f_{21}$ is everywhere regular and hence constant.

We now show:

\begin{prop}\label{exp-conns} The space of connections on $\E$ is given by
$f_{21}=C_1$, $f_{11}=c_1+c_2 x$, $f_{22}=c_3+c_4 x$,
and $f_{12}=c_5+c_6 x + c_7 x^2 +c_8 y +C_2 x^3$, where the
$c_i$ are arbitrary constants subject to the single linear 
relation $c_8 = C_2 (c_2-c_4)$, and $C_1 \text{ and } C_2$ are 
predetermined nonzero constants satisfying $C_1 C_2 =\frac{-1}{2}$.
\end{prop}

\begin{proof}We begin by looking at the lower right entry of the matrix 
for $\bar{T}_1$ in Equation \ref{exp-conneq},
and note the $\vp_{12}^{-1} \frac{d\vp_{12}}{\omega_1}$ has a simple
pole at $w$ which must be cancelled by one of the other terms. We
also note that since $\vp_\E=\vp_{12}^{-2}$, and $f_{21}$
must be constant, the term $\vp_{12} \vp_\E f_{21} =
\vp _{12}^{-1}f_{21}$ is regular on $U_1$ away from $w$, where it
can have at most a simple pole. Thus the $\vp_{12}^2 f_{22}$ term must
likewise be regular on $U_1$ away from $w$, with at most a simple pole
at $w$. Since $f_{22}$ must be regular on $U_2$ by hypothesis,
we conclude it is regular on $C$ except possibly for a pole of order at
most $3$ at $w$. But such a pole of order $3$ isn't
possible, so $f_{22} \in \Gamma(\O_C(2[w]))$. This means that the simple poles 
of the other two terms
must cancel, and $f_{21}$ is determined as a (nonzero) constant $C_1$:
explicitly, $C_1=\frac{d\vp _{12}}{\omega_1}(w)$. Precisely the same 
argument applies to the upper right entry, placing
$f_{11} \in \Gamma(\O_C(2[w]))$, so it only remains to analyze the upper right entry of the
matrix.

We immediately observe that on $U_1$, each term (excluding the $\vp_{12}^4
f_{12}$ term) is regular except possibly for a pole of order at most $2$
at $w$, which of course implies that $\vp_{12}^4 f_{12}$ is also, and
we can conclude that $f_{12}$ is regular on $C$ except for a pole of 
order at most $6$ at $w$. Then
we have $f_{21}=C_1 \in k^*$, $f_{11}=c_1+c_2 x$, $f_{22}=c_3+c_4 x$,
and $f_{12}=c_5+c_6 x + c_7 x^2 +c_8 y +C_2 x^3$, and we claim
that $C_2$ is also determined: the only other terms which
can have double poles are $-\vp_{12}^2\vp_\E ^2 f_{21}+\vp_\E 
\frac{d\vp_{12}}{\omega_1} - \vp_{12}\frac{d \vp_\E}{\omega_1}
=-\vp_{12}^{-2} f_{21}+3\vp_{12}^{-2}\frac{d\vp_{12}}{\omega_1}$
which are now predetermined, so $C_2$ is also determined, explicitly
as $-2 (\vp_{12}^{-6} x ^{-3})(w)C_1$. Lastly, 
we note that there is a linear relation on $c_2, c_4, \text{ and } 
c_8$ to insure that the simple poles cancel.

To conclude the proof, we use formal local analysis at $w$ to
obtain the desired statements on this linear relation and $C_1$ and $C_2$. 
Explicitly, our linear relation is given as 
$c_8=((\vp_{12}^{-3}y^{-1}x)(w))(c_2-c_4)+
((y^{-1}\vp_{12}^{-5})(w))((\vp_{12}^{-1}(C_1-3\frac{d\vp_{12}}{\omega_1})
-C_2 x^3 \vp_{12}^5)(w))$. Now, choose a local coordinate $z$ at $w$;
we will denote by $\ell_z(f)$ and $\ell_z'(f)$ the leading and second terms of the Laurent series expansion for $f$ in terms of $z$. From our relation between $x$ and $y$, we have $\ell_z(x)^5=\ell_z(y)^2$ and 
$2 \ell_z(y) \ell_z'(y) = 5 \ell_z(x)^4 \ell_z'(x)$. Simply considering leading terms, we find that since $\omega_1 = \vp_{12}^{-2}y^{-1}dx$, we have $C_1 = \frac{-\ell_z(\vp_{12})^{3}\ell_z(y)}{2\ell_z(x)}$, and $C_2 = 
\frac{\ell_z(x)}{\ell_z(\vp_{12})^3 \ell_z(y)}$. Thus, we have that
$C_1 C_2 = \frac{-1}{2}$, and also that $C_2$ is the coefficient of $(c_2-c_4)$ in our linear relation. It only remains to show that the constant term in that relation is in fact $0$. We may write it as 
$((y^{-1}\vp_{12}^{-5})(w))((\vp_{12}^{-1}(C_1-3\frac{d\vp_{12}}{\omega_1} 
-C_2 x^3 \vp_{12}^6))(w)$, so it suffices to show that 
$C_1- 3\frac{d\vp_{12}}{\omega_1} -C_2 x^3 \vp _{12}^6$, 
which we know must vanish at $w$, in fact vanishes to order
at least $2$ at $w$. For this, it is convenient to specialize to $z = \vp_{12}$, whereupon our earlier relation simplifies to 
$\ell_z(y)=C_2^{-1} \ell_z(x)$. We also compute $\ell_z(x)^3 = C_2^{-2}$, from which it follows that we can write $\ell_z'(y) = \frac{5}{2} C_2^{-1}
\ell_z'(x)$. We can now write everything in terms of
$\ell_z(x), \ell_z'(x)$ and $C_2$, and check directly that we get the desired cancellation to order $2$.
\end{proof}

We also consider the endomorphisms of $\E$, so that we can normalize
our connections via transport to simplify calculations. An endomorphism
is given by matrices $S_i$ regular on $U_i$, satisfying the relationship
$S_1 = E S_2 E^{-1}$. If we write $S_2 =
\begin{bmatrix} g_{11} & g_{12} \\ g_{21} & g_{22}
\end{bmatrix}$ we find that 
\begin{equation}\label{exp-end-eq}
S_1 = 
\begin{bmatrix}
g_{11} + \vp _{12} ^{-1} \vp_{\E} g_{21} & \vp_{12} ^2 g_{12} + 
\vp _{12} \vp _{\E} g_{22} - \vp _{12} \vp _{\E} g_{11} - \vp_{\E}^2 g_{21} \\
\vp _{12} ^{-2} g_{21} & g_{22} - \vp_{12} ^{-1} \vp _{\E} g_{21}
\end{bmatrix}
\end{equation}

We can now compute directly:

\begin{prop}\label{exp-ends} The space of endomorphisms of $\E$ is given by
$g_{21}=0$, $g_{11}=g_{22} \in k$, and 
$g_{12} = g_{12}^0 + g_{12}^1 x \in \Gamma(\O_C(2[w]))$.
Every connection on $\E$ has a unique transport-equivalent connection
with $f_{11}=0$.
\end{prop}

\begin{proof}Noting that the lower left entry for $S_1$ 
in equation \ref{exp-end-eq} is $\vp _{12} ^{-2} g_{21}$, we see that
$g_{21}$ has to be regular everywhere on $C$, and vanishes to order
at least $2$ at $w$; hence, it is $0$. We then see that the upper
left and lower right entries are just $g_{11}$ and $g_{22}$
respectively, meaning that these are both everywhere regular and
hence constant. Finally, the upper right term is then 
$\vp _{12} ^2 g_{12} + \vp _{12} ^{-1} (g_{22}-g_{11})$; the second
term will have a simple pole at $w$ if and only if $g_{22} \neq g_{11}$,
and since $g_{12}$ cannot have a triple pole at $w$, we conclude that
$g_{22} = g_{11}$, and finally that $g_{12} \in \Gamma(\O_C(2[w]))$, giving the description of the endomorphisms of $\E$.

Such an endomorphism is invertible if and only if $g_{11} \neq 0$. Since
transport along an automorphism is invariant under scaling the automorphism, we can then set $g_{11}=g_{22}=1$ without loss of generality.
Now, since $S_2$ is upper triangular, with constant diagonal coefficients,
$S_2^{-1} \frac{dS_2}{\omega_2}$ has only its upper right coefficient
non-zero. Moreover, conjugating $\bar{T}_2$ by $S_2$ will simply substract
$f_{21} g_{12}$ from the upper left coefficient of $\bar{T}_2$. Since we
know $f_{21}$ is a determined nonzero constant, and $g_{12}$ and $f_{11}$
can both be arbitrary in $\Gamma(\O_C(2[w]))$, this means that each 
connection has a unique transport class with $f_{11}=0$, as desired.
\end{proof}

Thus, from now on we will normalize our calculations as follows:
set $f_{11} = 0$ by transport; set $f_{22} =0$ since we want the 
determinant connection (obtained by taking the trace) to be $0$; and set $f_{21} =1$. We accomplish the last by scaling $\vp_{12}$ appropriately: we 
saw that $f_{21}= \frac{d\vp_{12}}{\omega_1}(w)$, and recalling that 
$\omega_1 = \vp_{12}^{-2} y^{-1} dx$, it suffices to scale $\vp_{12}$ by a
cube root of $f_{21}$. We also note that this does 
not pose any problems for our prior choice of $\vp_{12} = \frac{x^2}{y}$;
one can check that for this choice, we have $f_{21}$ invariant as 
$\frac{-1}{2}$, and the scaling factor for $\vp_{12}$ is
independent of the $a_i$. Lastly, since $c_8=0$ now that $c_2=c_4=0$,
we conclude that we are reduced to considering the case:

\begin{sit}\label{exp-normalized} Our connection matrix $T_2$ on $U_2$ is of the form $T_2 =
\begin{bmatrix}0 & f_{12} \\ 1 & 0 \end{bmatrix}$, with
$f_{12} = c_5+c_6 x + c_7 x^2 - \frac{1}{2}x^3$.
\end{sit}

Finally, for later use we formally generalize our results.

\begin{prop}\label{exp-gen} Propositions \ref{exp-conns} and \ref{exp-ends} hold in the following more more general settings:
\begin{ilist}
\itm After base change to an arbitrary $k$-algebra $A$, if we replace the $k$-valued constants by $A$-valued constants;
\itm When we allow our defining polynomial $g(x)$ to degenerate to produce nodes away from $w$, if we replace $\Omega^1_C$ by the dualizing sheaf $\omega_C$ in the definition of connections;
\itm When we consider families of curves obtained from maps 
$k[a_1, \dots a_5] \rightarrow A$ taking values in the open subset $U_{\nod} \subset \A^5$ corresponding to at worst nodal curves.
\end{ilist}
\end{prop}

\begin{proof} For (i), if we denote by $f$ the map $\Spec A
\rightarrow \Spec k$, and $\pi$ the structure map $C \rightarrow \Spec k$, 
this coefficient replacement corresponds to the natural map $f^* \pi_* \F \rightarrow
\pi_{f*} f_{\pi}^* \F$ for the sheaves $\cEnd(\E) \otimes \Omega^1_C$ and
$\cEnd(\E)$. But since the base is a point, every base change is flat, and
it immediately follows \cite[Prop. III.9.3]{ha1} that this natural map is 
always an isomorphism, giving the desired statement.

For (ii), we need only note that our arguments go through unmodified, since $\omega_C$ is still isomorphic to $\O(2[w])$, and the same standard Riemann-Roch argument as in the smooth case still shows it that there can 
be no function in $\Gamma(\O(3[w]))\smallsetminus \Gamma(\O(2[w]))$.

Finally, for (iii) we make use of the fact that, as remarked
immediately above, we can choose $\vp_{12}$ to be a specific function 
varying algebraically in the whole family. Once again, if we denote by $\F$ the sheaf 
$\cEnd(\E) \otimes \omega_C$ or $\cEnd(\E)$ as appropriate, but this time 
in the universal setting over $U_{\nod}$, the theory of cohomology and base 
change gives that since $h^0(C, \F)$ is constant on fibers, $\pi_* \F$ is 
locally free of the same rank, and pushforward commutes with base change. 
Now, if we let our constants describing sections of $\F$ lie in
$k[a_1, \dots, a_5]$, we clearly obtain a subsheaf of $\pi_* \F$ of the correct
rank; further, the inclusion map is an isomorphism when restricted to every
fiber, so it must in fact be an isomorphism, which yields the desired 
result for arbitrary $A$ via base change.
\end{proof}

It follows formally that the closed subschemes we describe explicitly
corresponding to vanishing $p$-curvature in Section \ref{s-exp-pcurve} and 
nilpotent $p$-curvature in Section \ref{s-exp-det} are also functorial 
descriptions which hold for nodal curves.

\section{Calculations of $p$-curvature}\label{s-exp-pcurve}

Continuing with the situation and notations of the previous section, and in
particular that of Situation \ref{exp-presented}, we
conclude with the $p$-curvature calculations to complete the proof of
Theorem \ref{exp-main} for $p \leq 7$, except for the statement on the general curve in
characteristic $7$, which depends on the results of the subsequent section.

We write:
$$\psi_\nabla (\theta) = 
\begin{bmatrix}
h_{11} &
h_{12} \\
h_{21} &
h_{22} 
\end{bmatrix}$$

The first case we handle is $p=3$. Equation \ref{exp-ftheta-3} gave us
$f_{\theta ^3} = x^3 +a_3$. 
We show:

\begin{prop}In characteristic $3$, $\E$ has a unique transport 
equivalence class of 
connections with $p$-curvature zero and trivial determinant.
\end{prop}

\begin{proof}With all of our normalizations from Situation
\ref{exp-normalized}, the $p$-curvature matrix given by Equation
\ref{exp-pcurve-3} becomes rather tame:

$$\psi_\nabla (\theta) = 
\begin{bmatrix}
\theta f_{12} &
f_{12}^2 + \theta ^2 f_{12} - f_{\theta^3} f_{12} \\
f_{12}- f_{\theta^3} &
-\theta  f_{12} 
\end{bmatrix}$$

Even better, we note that we have 
$$h_{12} = \theta (h_{11}) + f_{12} h_{21},$$ 
so $h_{12}$ vanishes if $h_{11}$ and $h_{21}$ do. Similarly, recalling 
that by Lemma \ref{exp-thetaftheta}, $\theta f_{\theta ^3} = 0$, we see 
that $h_{11} = \theta (h_{21})$, and $h_{22} = - h_{11}$. Hence, to check 
if the $p$-curvature vanishes, it suffices to check that $h_{21}$ vanishes.

But this is a triviality, as we simply get that $h_{21}=0$ if and only if 
$f_{12} = a_3 + x^3$. Recalling that after normalization $f_{12}$ was 
given by $c_5 + c_6 x + c_7 x^2 - \frac{1}{2} x^3$, we get the unique solution $c_5 = a_3, c_6=c_7=0$.
\end{proof}

We now handle the case $p=5$. We had from Equation \ref{exp-ftheta-5} that
$f_{\theta^5} = 2 a_1 x^5 + a_3^2 + 2 a_2 a_4 + 2 a_1 a_5$.

\begin{prop}In characteristic $5$, the number of transport equivalence
classes of connections with $p$-curvature zero and trivial determinant
is given as the number of roots of the quintic polynomial:
\begin{multline*}
(3 a_1 a_2 ^2 + 3 a_2 a_3 + a_5) + 
(a_1 ^2 a_2 + a_2 ^2 + 3 a_1 a_3 + 4 a_4) c_5 \\
+ (3 a_1 ^3 + 4 a_1 a_2 + a_3) c_5 ^2 +
(3 a_1 ^2 + 4 a_2) c_5 ^3 + a_1 c_5 ^4 + 4 c_5 ^5
\end{multline*}
\end{prop}

\begin{proof}
With our normalizations as above, in terms of $f_{12}$ and $f_{\theta^5}$, 
the $p$-curvature matrix obtained from Equation \ref{exp-pcurve-5} is 

$$\psi_\nabla (\theta) = 
\begin{bmatrix}
4 f_{12} \theta (f_{12}) + \theta^3(f_{12})&
f_{12}^3 + 4 (\theta(f_{12}))^2 + 2 f_{12} \theta^2(f_{12}) 
+ \theta^4 (f_{12}) + 4 f_{12} f_{\theta^5} \\
f_{12}^2 + 3 \theta^2(f_{12}) + 4 f_{\theta ^5} &
f_{12} \theta(f_{12}) + 4 \theta^3 (f_{12})
\end{bmatrix}$$

Conveniently, we note that as before it actually suffices to check that 
$h_{21}$ is $0$, since we see that $h_{22}= 3 \theta _0 (h_{21})$, that 
$h_{11} = -h_{22}$, and that 
$h_{12} = f_{12} h_{21} + 2 \theta ^2 (h_{21})$.

Substituting in for $f_{12}$ and $f_{\theta^5}$, we get that the 
remaining (lower left) term is given by

\begin{multline*}
(4 a_3 ^2 + 3 a_2 a_4 + 3 a_1 a_5 + c_3 ^2 + a_5 c_5) + 
(a_5 + 3 a_3 c_4 + 2 c_3 c_4 + 4 a_4 c_5) x \\
+ (2 a_2 c_4 + c_4 ^2 + 2 a_3 c_5 + 2 c_3 c_5) x ^2 \\
+ (4 a_3 + 4 c_3 + a_1 c_4 + 2 c_4 c_5) x ^3 
+ (3 a_2 + 4 c_4 + 3 a_1 c_5 + c_5^2) x ^4
\end{multline*}

Setting the $x^3$ and $x^4$ terms to $0$ allows us to solve for $c_4$ and $c_3$. Substituting in, we find that the $x^2$ term drops
out, while the coefficient of $x$ is:

\begin{multline*}
(3 a_1 a_2 ^2 + 3 a_2 a_3 + a_5) + 
(a_1 ^2 a_2 + a_2 ^2 + 3 a_1 a_3 + 4 a_4) c_5 \\
+ (3 a_1 ^3 + 4 a_1 a_2 + a_3) c_5 ^2 +
(3 a_1 ^2 + 4 a_2) c_5 ^3 + a_1 c_5 ^4 + 4 c_5 ^5
\end{multline*}

The constant coefficient is $c_5+ 3 a_1$ times the $x$ coefficient,
so we get that the connections with $p$-curvature $0$ correspond
precisely to the roots of the above polynomial, as asserted.
\end{proof}

Lastly, we take a look at the case $p=7$. Equation \ref{exp-ftheta-7} gave
us:
$$f_{\theta^7} = a_3 ^3 + 6 a_2 a_3 a_4 + 3 a_1 a_4 ^2 + 3 a_2 ^2 a_5 + 6
a_1 a_3 a_5 +
6 a_4 a_5 + (3 a_1 ^2 + 3 a_2) x ^7.$$

We will show: 

\begin{prop}In characteristic 7, the number of transport equivalence classes
of connections on $\E$ with $p$-curvature $0$ and trivial determinant is 
given as the intersection of four plane curves in $\A^2$. For a general
curve, it is positive. The locus $F_{2,7}$ of transport equivalence classes 
of connections on $\E$ with $p$-curvature $0$ and trivial determinant 
considered over the $\A^5$ with which we parametrize genus $2$ curves
is cut out by $4$ hypersurfaces in $\A^5 \times \A^2$.
\end{prop}

\begin{proof}
Here, even with our normalizations the $p$-curvature matrix obtained from
Equation \ref{exp-pcurve-7} is rather messy, but we find its coefficients 
are given by:

$$h_{11} = 2 f_{12} ^2 \theta (f_{12}) + \theta (f_{12}) \theta ^2 (f_{12}) -
3 f_{12} \theta ^3 (f_{12}) + \theta ^5 (f_{12})$$

$$h_{21} = -f_{\theta^7} + f_{12} ^3 + 3 (\theta (f_{12})) ^2 - f_{12} \theta ^2
(f_{12}) - 2 \theta ^4 (f_{12})$$

$$h_{12} = -f_{\theta^7} f_{12} + f_{12} ^4 + f_{12} ^2 \theta ^2 (f_{12}) +
(\theta ^2 (f_{12})) ^2 -
2 \theta (f_{12}) \theta ^3 (f_{12}) + 2 f_{12} \theta ^4 (f_{12}) + \theta ^6
(f_{12})$$

$$h_{22} = -2 f_{12} ^2 \theta (f_{12}) - \theta (f_{12}) \theta ^2 (f_{12}) 
+ 3 f_{12} \theta ^3 (f_{12}) - \theta ^5 f_{12}$$

Once again, it is enough to consider a single one of these coefficients,
as we see that $h_{11} = 3 \theta(h_{21})$, that 
$h_{12} = f_{12} h_{21} + 3 \theta^2(h_{21})$, and that
$h_{22} = -h_{11}$.

Looking then at the formula for $h_{21}$, substituting in for $f_{12}$ 
and $f_{\theta^7}$ gives a polynomial of
degree $6$ in $x$. The $x^6$ term lets us
solve for $c_3$:

$$c_3 = 5 a_1 a_2 +
a_3 + 4 a_1 c_4 + 4 a_1^2 c_5 + c_4 c_5 + 2 a_1 c_5^2 + 5 c_5^3$$

The $x^5$ term is then

\begin{multline}h_{7,1}=2 a_1 ^2 a_2 + a_1 a_3 + 5 a_4 + 4 a_1 ^2 c_4 + 
5 a_2 c_4 + 6 c_4 ^2 \\
+ 3 a_1 ^3 c_5 + 6 a_1 a_2 c_5 + 3 a_3 c_5 + 5 a_1 c_4 c_5 + 3 a_1 c_5 ^3 
+ 6 c_5 ^4.
\end{multline}
while the $x^4$ term is $-c_5$ times the $x^5$ term, and the $x^3$ term is
$-(c_5^2+a_1 c_5 + 3a_2 + c_4)$ times the $x^5$ term. Taking the $x^2$ term
minus $-(5 c_5^3+ 5 a_1 c_5^2 + 2 c_4 c_5 + 5 a_1 a_2 + 4 a_3 + 2 a_1 c_4)$
times
the $x^5$ term leaves:

\begin{multline}h_{7,2}= 3 a_1 ^3 a_2 ^2 + 6 a_1 ^2 a_2 a_3 + 4 a_1 a_3 ^2 
+ 4 a_3 a_4 + 2 a_2 a_5 + 3 a_1 ^3 a_2 c_4 + 4 a_1 ^2 a_3 c_4 \\
+ 2 a_1 a_4 c_4 + 4 a_5 c_4 + a_1 ^3 c_4 ^2 + a_3 c_4^2 + 3 a_1 c_4^3 
+ a_1 ^4 a_2 c_5 \\
+ 5 a_1 ^3 a_3 c_5 + a_1 ^2 a_4 c_5 + 3 a_1 a_5 c_5 + 6 a_1 ^4 c_4 c_5 
+ a_1 ^2 a_2 c_4 c_5 \\
+ a_1 a_3 c_4 c_5 + 3 a_4 c_4 c_5 + a_1 ^2 c_4 ^2 c_5 + 5 a_2 c_4 ^2 c_5 
+ c_4 ^3 c_5 \\
+ 4 a_1 ^3 a_2 c_5 ^2 + 6 a_1 ^2 a_3 c_5 ^2 + a_1 a_4 c_5 ^2 + 3 a_5 c_5 ^2 
+ 3 a_1 ^3 c_4 c_5 ^2 \\
+ a_1 a_2 c_4 c_5 ^2 + a_3 c_4 c_5 ^2 + 3 a_1 ^2 a_2 c_5 ^3 
+ a_1 a_3 c_5 ^3 + 4 a_1 ^2 c_4 c_5 ^3 + 6 c_4 ^2 c_5 ^3.
\end{multline}

Similarly, taking the $x$ term minus $-(5 c_4 c_5^2+5 a_1 c_4 c_5 +6 a_4 + 2
c_4^2)$
times the $x^5$ term leaves:

\begin{multline}h_{7,3} = 5 a_1 ^2 a_2 a_4 + 6 a_1 a_3 a_4 + a_1 a_2 a_5 
+ 5 a_1 ^2 a_2 ^2 c_4 + 4 a_1 a_2 a_3 c_4 + 3 a_1 ^2 a_4 c_4 \\
+ 2 a_1 a_5 c_4 + 5 a_1 ^2 a_2 c_4 ^2 + a_1 a_3 c_4 ^2 + 3 a_4 c_4 ^2 
+ 3 a_2 c_4 ^3 \\
+ 5 c_4 ^4 + 4 a_1 ^3 a_4 c_5 + 6 a_3 a_4 c_5 + 5 a_1 ^2 a_5 c_5
+ 3 a_2 a_5 c_5 \\
+ 4 a_1 ^3 a_2 c_4 c_5 + 4 a_1 ^2 a_3 c_4 c_5 + a_1 a_4 c_4 c_5 
+ 6 a_5 c_4 c_5 \\
+ 3 a_1 ^3 c_4 ^2 c_5 + 4 a_1 a_2 c_4 ^2 c_5 + 4 a_3 c_4 ^2 c_5 
+ a_1 c_4 ^3 c_5 \\
+ 2 a_1 ^2 a_4 c_5 ^2 + 6 a_1 a_5 c_5 ^2 + 2 a_1 ^2 a_2 c_4 c_5 ^2 
+ 2 a_1 a_3 c_4 c_5 ^2 \\
+ a_4 c_4 c_5 ^2 + 5 a_1 ^2 c_4 ^2 c_5 ^2 + 4 a_2 c_4 ^2 c_5 ^2 
+ 5 c_4 ^3 c_5 ^2 \\
+ 5 a_1 a_4 c_5 ^3 + a_5 c_5 ^3 + 5 a_1 a_2 c_4 c_5 ^3 + 5 a_3 c_4 c_5 ^3 
+ 2 a_1 c_4 ^2 c_5 ^3.
\end{multline}

Lastly, taking the constant term minus 
\begin{multline*}-(6 c_5^5+5 c_4 c_5^3 + 3 a_1 ^2 c_5^3 +2 a_1^3 c_5^2 
+ 5 a_1 a_2 c_5^2 + 2 a_3 c_5^2 + 6 a_1 c_4 c_5^2 \\
+ 5 a_1^2 a_2 c_5 + 2 a_1 a_3 c_5 + 2 a_4 c_5 + a_1^2 c_4 c_5 \\
+ 2 a_2 c_4 c_5 + 2 c_4^2 c_5 + 6 a_1 a_4 
+ 4 a_5 + 4 a_1 a_2 c_4 + 3 a_3 c_4 + 3 a_1 c_4^2)
\end{multline*}
times the $x^5$ term leaves:

\begin{multline}h_{7,4}=6 a_1 ^3 a_2 ^3 + 5 a_1 ^2 a_2 ^2 a_3 
+ a_1 a_2 a_3 ^2 + 5 a_1 ^3 a_2 a_4 + 6 a_1 ^2 a_3 a_4 + a_2 a_3 a_4 \\
+ 6 a_1 a_4 ^2 + a_1 ^2 a_2 a_5 + 4 a_2 ^2 a_5 + 5 a_1 a_3 a_5 + 4 a_4 a_5\\ 
+ 4 a_1 ^2 a_2 a_3 c_4 + a_1 a_3 ^2 c_4 + 3 a_1 ^3 a_4 c_4 
+ 2 a_1 a_2 a_4 c_4 \\
+ 3 a_3 a_4 c_4 + 2 a_1 ^2 a_5 c_4 + 3 a_1 ^3 a_2 c_4 ^2 
+ 6 a_1 a_2 ^2 c_4 ^2 \\
+ a_2 a_3 c_4 ^2 + 6 a_5 c_4 ^2 + 6 a_1 ^3 c_4 ^3 + 4 a_1 a_2 c_4 ^3 
+ 4 a_3 c_4 ^3 \\
+ 4 a_1 c_4 ^4 + 2 a_1 ^4 a_2 ^2 c_5 + 3 a_1 ^3 a_2 a_3 c_5 
+ 4 a_1 ^4 a_4 c_5 \\
+ 2 a_1 ^2 a_2 a_4 c_5 + 2 a_1 a_3 a_4 c_5 + 5 a_1 ^3 a_5 c_5 + a_3 a_5 c_5 
+ 3 a_1 ^4 a_2 c_4 c_5 \\
+ 2 a_1 ^2 a_2 ^2 c_4 c_5 + 2 a_1 ^3 a_3 c_4 c_5 + 2 a_1 a_2 a_3 c_4 c_5 
+ 5 a_3 ^2 c_4 c_5 \\
+ 6 a_1 ^2 a_4 c_4 c_5 + 6 a_2 a_4 c_4 c_5 + 5 a_1 a_5 c_4 c_5 
+ 2 a_1 ^4 c_4 ^2 c_5 + 2 a_1 ^2 a_2 c_4 ^2 c_5 \\
+ 3 a_2 ^2 c_4 ^2 c_5 + 6 a_1 a_3 c_4 ^2 c_5 + 4 a_4 c_4 ^2 c_5 
+ a_2 c_4 ^3 c_5 + 5 c_4 ^4 c_5 \\
+ a_1 ^3 a_2 ^2 c_5 ^2 + 5 a_1 ^2 a_2 a_3 c_5 ^2 + 2 a_1 ^3 a_4 c_5 ^2 
+ 2 a_1 a_2 a_4 c_5 ^2 \\
+ 2 a_3 a_4 c_5 ^2 + 6 a_1 ^2 a_5 c_5 ^2 + 5 a_1 ^3 a_2 c_4 c_5 ^2 
+ 2 a_1 a_2 ^2 c_4 c_5 ^2 \\
+ a_1 ^2 a_3 c_4 c_5 ^2 + 2 a_2 a_3 c_4 c_5 ^2 + 4 a_1 a_4 c_4 c_5 ^2
+ 5 a_5 c_4 c_5 ^2 \\
+ a_1 ^3 c_4 ^2 c_5 ^2 + 6 a_1 a_2 c_4 ^2 c_5 ^2 + 2 a_1 c_4 ^3 c_5 ^2 
+ 6 a_1 ^2 a_2 ^2 c_5 ^3 \\
+ 2 a_1 a_2 a_3 c_5 ^3 + 5 a_1 ^2 a_4 c_5 ^3 + a_1 a_5 c_5 ^3 
+ 2 a_1 ^2 a_2 c_4 c_5 ^3 \\
+ 6 a_1 a_3 c_4 c_5 ^3 + 5 a_4 c_4 c_5 ^3 + 6 a_1 ^2 c_4 ^2 c_5 ^3 
+ 4 a_2 c_4 ^2 c_5 ^3 + 3 c_4 ^3 c_5 ^3.
\end{multline}

These four polynomials are then the defining equations in characteristic
$7$, describing the locus as an intersection of $4$ affine plane curves,
as desired. By direct computation in Macaulay 2, the coordinate ring of the 
affine algebraic set cut out by these equations has dimension 5. Since we know that it can only have dimension $0$ over any given choice for the 
$a_i$, this implies that it has a non-empty fiber for a general choice of
$a_i$, yielding the positivity assertion.
\end{proof}

Finally, we compute an example which will allow us to deduce the characteristic $7$ case of Theorem \ref{exp-main} in the next section.

\begin{lem}\label{exp-ch7-ex} For the curve given by $a_1=a_2=a_3=0, a_4=1$, and $a_5=3$, there
are $14$ solutions to our equations, all reduced. Further, the local rings
of $F_{2,7}$ at each of these points are all isomorphic. 
\end{lem}

\begin{proof}
First, we set $a_1=a_2=a_3=0$, $a_4=1$ and $a_5=3$, and our defining equations
become considerably simpler:

$$h_{7,1}= 5 + 6 c_4^2 + 6 c_5^4$$
$$h_{7,2}= 5 c_4 + 3 c_4 c_5 + c_4^3 c_5 + 2 c_5^2 + 6 c_4 ^2 c_5^3$$
$$h_{7,3}= 3 c_4^2 + 5 c_4^4 + 4 c_4 c_5 + c_4 c_5^2 + 5 c_4^3 c_5^2 + 3 c_5^3$$
$$h_{7,4}= 5+ 4 c_4^2 + 4 c_4 ^2 c_5 + 5 c_4 ^4 c_5 + c_4 c_5^2 + 5 c_4 c_5^3 
+ 3 c_4^3 c_5^3$$

If we use $h_{7,1}$ to substitute for $c_4^2$ in $h_{7,2}$, we get:

$$c_4(5+c_5+ 6 c_5^5) + c_5^2(2+ 2 c_5 + c_5^5)$$

We check that we cannot have $5 + c_5 + 6 c_5^5=0$, so we can localize away from $5+c_5 + 6 c_5^5$, 
setting $c_4 = \frac{c_5^2(2+ 2 c_5 + c_5^5)}{5 + c_5 + 6 c_5^5}$. Making
this substitution and taking numerators, the $h_{7,i}$ give four
polynomials in $c_5$. However, they are multiples of the polynomial
given by $h_{7,1}$, which is:

$$6 + c_5 + 5 c_5 ^2 + 6 c_5 ^4 + 2 c_5 ^5 + 6 c_5 ^6 + 6 c_5 ^9 + 3 c_5
^{10} + 5 c_5 ^{14}$$

This then gives the $14$ reduced solutions, and the fact that the local
rings of $F_{2,7}$ at each of these points are isomorphic follows from
the fact that this degree $14$ polynomial is irreducible over $\fF_7$, since
the $14$ points are then Galois conjugate in $F_{2,7}$, which is defined over
$\fF_7$.
\end{proof}

\section{On The Determinant of the $p$-Curvature Map}\label{s-exp-det}

In this section we explicitly calculate the highest degree terms of 
$\det \psi$, 
the determinant of the $p$-curvature map in the case of a genus 2 curve and
the specific unstable vector bundle of Situation \ref{exp-specific-e}. We 
use the calculation to prove that $\det \psi$ is 
finite flat, of degree $p^3$, and therefore conclude that in families of 
curves, the kernel of $\det \psi$ is finite flat. This has immediate 
implications for the connections on $\E$ of vanishing $p$-curvature as well, in particular allowing us to finish the proof of the characteristic-specific portion of Theorem \ref{exp-main}. The results here are a special case of Mochizuki's work (see \cite[Thm. II.2.3, p. 129]{mo1}), obtained by an argument which is essentially the same, but discovered independently, and significantly simpler in the special case handled here.

We wish to compute in our specific situation the morphism $\det \psi^0$ 
(Proposition \ref{exp-pcurve-formal} (iv)), which is to say, the morphism obtained from $\psi^0$
(Proposition \ref{exp-pcurve-formal} (iii)) by taking the determinant. In fact, we take $\psi^0$ 
to be the induced map on transport-equivalence classes of connections.
We remark that in the situation of rank $2$ vector bundles with 
trivial determinant, and after restricting to connections with trivial 
determinant, because the image of $\psi^0$ is contained among the traceless
endomorphisms, the vanishing of the determinant of the
$p$-curvature is then equivalent to nilpotence of the endomorphisms given 
by the $p$-curvature map. Such connections are frequently called
{\bf nilpotent} in the literature (see, for instance, \cite{ka1} or
\cite{mo3}).

We now take our curve $C$ of genus $2$ from before, with $\E$ the particular
unstable bundle of rank $2$ we had been studying, as in
Situations \ref{exp-presented}, \ref{exp-specific-e}, and
\ref{exp-normalized}. We also take the 
particular $\theta$ from before, with $\hat{\theta}(\omega_2)=1$. Since 
$\omega_2$ has a double zero at $w$, we see that $\theta$ has a double pole 
there, so that our explicit identification of $\Omega^1_C$ is as 
$\O(2 [w])$. We know that our space of connections with trivial determinant 
on $\E$ is (modulo transport) $3$-dimensional, and of course 
$h^0(C^{(p)}, (\Omega^1_{C^{(p)}})^{\otimes 2})=
\deg (\Omega^1_{C^{(p)}})^{\otimes 2}+1-g= 4g-4+1-g = 3g-3 = 3$, so we have a map 
from $\A^3$ to $\A^3$. We choose coordinates
on the first space to be given by the $(c_5, c_6, c_7)$ determining $f_{12}$, 
while 
the function we will get will be of the form $f_1(c_5, c_6, c_7)+f_2(c_5, c_6,
c_7)x^p+f_3(c_5, c_6, c_7)x^{2p}$, and we obtain coordinates on the image 
space as the monomials $(1,x^p, x^{2p})$.

We will use our earlier calculations to recover, in a completely explicit 
and elementary fashion, the genus $2$ case of Mochizuki's result:

\begin{thm}\label{exp-nilp} On the unstable vector bundle $\E$ described by Situation
\ref{exp-specific-e} for a smooth proper genus $2$ curve $C$ as in Situation
\ref{exp-presented}, the map $\det \psi^0$ is a finite flat morphism from 
$\A^3$ to $\A^3$, of degree $p^3$. Further, $\det \psi^0$ remains finite flat
when considered as a family of maps over the open subset $U_{\ns} \subset 
\A^5$ corresponding to nonsingular curves. Lastly, the induced
map from $\ker \det \psi^0$ to $U_{\ns}$ is finite flat.
\end{thm}

\begin{proof} It suffices to prove the asserted finite flatness for the
family of maps $\A^3 \times U_{\ns} \rightarrow \A^3 \times U_{\ns}$ over
$U_{\ns}$,
since the statements on individual curves and on the kernel of $\det
\psi$ both follow from restriction to fibers. This is turn will follow from 
the claim that the leading term of $f_i$ is $-c_{i+4}^p$, with 
all other terms of strictly lesser total degree in the $c_i$. We prove this by direct calculation.

If $T = 
\begin{bmatrix}
0 & f_{12} \\
1 & 0
\end{bmatrix}$
is the connection matrix for $\nabla$, we claim that the leading term
will come from the $T^p$ term in the $p$-curvature formula. Now, $T^2 =
\begin{bmatrix}
f_{12} & 0\\
0 & f_{12}
\end{bmatrix}$, so we find 
$$T^p = 
\begin{bmatrix}
0 & (f_{12})^{\frac{p+1}{2}} \\
(f_{12})^{\frac{p-1}{2}} & 0
\end{bmatrix}$$

Next, $f_{12}$ is linear in the $c_i$, as are $\theta ^i f_{12}$ for all $i$.
Considering the $p$-curvature formula coefficients as polynomials in
$\theta ^i f_{12}$, we will show that the degree of the remaining terms
are all less than or equal to $\frac{p-1}{2}$, with the degree of the
terms in the lower left strictly less. This will imply that the leading
term of the determinant is given by 
$$-(f_{12})^p = -(c_5+c_6 x + c_7 x^2 - \frac{1}{2}x^3)^p =
-c_5^p - c_6^p x^p - c_7^p x^{2p}+ \frac{1}{2^p}x^{3p}$$ 
giving the desired formula for the leading terms of the constant, $x^p$, 
and $x^{2p}$ terms.

We observe that since $\theta^i T =
\begin{bmatrix}
0 & \theta ^i f_{12} \\
0 & 0
\end{bmatrix}$
for all $i>0$, $(\theta^i T)( \theta ^j T)=0$ for any $i,j>0$. We use this
and the fact that $T^2$ is diagonal to write any term in the $p$-curvature
as one of the following:
\begin{nlist}
\itm $T^{2 i_0} (\theta^{i_1}T) T \dots (\theta^{i_k}T)$ 
\itm $T^{2 i_0} T (\theta^{i_1}T) T \dots (\theta^{i_k}T)$ 
\itm $T^{2 i_0} (\theta^{i_1}T) T \dots (\theta^{i_k}T) T$ 
\itm $T^{2 i_0} T (\theta^{i_1}T) T \dots (\theta^{i_k}T) T$ 
\end{nlist}
where $2 i_0 + \sum _{j>0} (i_j+2) = p+1, p, p, p-1$ respectively.

We observe that these correspond to non-zero upper right, lower right,
upper left, and lower left coefficients, respectively (in particular,
at most one is non-zero). We know that the first term is a scalar matrix of 
degree $i_0$ in $f_{12}$. We see that $T (\theta^{i_j} T) = 
\begin{bmatrix}
0 & 0 \\
0 & \theta ^i f_{12} 
\end{bmatrix}$,
so a product of $k-1$ such terms has total degree $k-1$ in the
$\theta ^i f_{12}$. Lastly, multiplying on the left by $(\theta ^{i_1}T)$
raises the degree by one and moves the nonzero coefficient back to the
upper right. Thus, in the first case, we get total degree $i_0+k$. But we
see that this is actually the same in the other cases, as multiplying on
the left or right by $T$ just moves the nonzero coefficient, without 
changing it. Finally, with $k>0$, we have $i_0 +k < \frac{1}{2}( 2 i_0 +
\sum _{j>0} (i_j+2)) $, which is $\frac{1}{2}$ times $p+1, p, p$ or $p-1$ 
depending
on the case. But this is precisely what we wanted to show, since it
forces the degree to be less than or equal to $\frac{p-1}{2}$ in the 
first three cases, and strictly less in the fourth.

Lastly, $-f_{\theta^p} T$ is linear in the $c_i$ in the upper right term,
and constant in the rest, so doesn't cause any problems for $p\geq 3$.
\end{proof}

We can immediately conclude:

\begin{cor}\label{exp-finite} The subscheme of $U_{\ns} \times \A^3$ giving 
connections with $p$-curvature $0$ is finite over $U_{\ns}$.
\end{cor}

We are now ready to put together previous results to finish the proof of our
main theorem in the case of characteristic $7$:

\begin{proof}[Proof of Theorem \ref{exp-main}, $p=7$ case] We simply apply our finiteness result to our explicit example from Lemma \ref{exp-ch7-ex}.
We calculated that $F_{2,7}$ has dimension $5$, so by properness the local 
ring of at least one point in our example has dimension $5$, hence they all 
do. By the reducedness of our example, all its points are unramified over 
the base, and by finiteness, we conclude that on an open subset of the base containing our chosen point, $F_{2,7}$ is finite and unramified, and everywhere $5$-dimensional. Then, by the regularity of the base, we find that over this open set, $F_{2,7}$ must be regular, hence flat, hence \'etale, so we conclude the desired statement for a general curve in characteristic $7$.
\end{proof}

\section{Connections and Nodes}\label{s-def-background}

This section and the next draw heavily on the results and ideas of Sections
2 and 3 of \cite{os10}.

In this section, we discuss connections on nodal curves, and classify them
in terms of gluings of connections on the normalization. For the sake of
simplicity and generality, we follow Mochizuki's argument for the gluing,
with the only difference being that because we are not working with
projective bundles, we must rigidify our situation by specifying glued line sub-bundles $\L$, as in Proposition \ref{def-glue-prop}.

Let $C$ be a proper nodal curve, and $\E$ a vector bundle on $C$. We begin by fixing some terminology:

\begin{defn}A {\bf logarithmic connection} on $\E$ is a $k$-linear map 
$\nabla: \E
\rightarrow \E \otimes \omega _C$, where $\omega_C$ is the dualizing sheaf
on $C$, and $\nabla$ satisfies the Liebnitz rule induced by the canonical 
map $\Omega^1_C \rightarrow \omega_C$. Given a reduced divisor $D$ supported on the smooth locus of $C$, a {\bf $D$-logarithmic connection} on $\E$ is a $k$-linear map $\nabla: \E \rightarrow \E \otimes \omega _C(D)$ satisfying the Liebnitz rule.
\end{defn}

We note that with the exception of the Cartier isomorphism, all the properties of connections and $p$-curvature which we have used
still hold if one replaces $\Omega^1_C$ by $\omega_C$ (and in
particular, the sheaf of derivations by $\omega_C^\vee$) throughout. We
summarize as follows.

\begin{prop}\label{def-background} All statements on induced connections for
operations of vector bundles, and the basic properties of the $p$-curvature map of Proposition \ref{exp-pcurve-formal}, hold in the case of 
logarithmic connections on nodal curves, with $\omega_C$ in place of 
$\Omega_C^1$.
One still has a canonical connection on a Frobenius pullback with 
vanishing $p$-curvature whose kernel recovers the original sheaf.
\end{prop}

Although it is true that taking the kernel of the canonical connection
of a Frobenius pullback still recovers the original sheaf on $C^{(p)}$ when
$C$ is singular, the Cartier isomorphism fails because given a logarithmic
connection with vanishing $p$-curvature on $C$, the Frobenius pullback of
the kernel will not in general map surjectively onto the original sheaf at 
the singularities of $C$.

\begin{notn}Let $\tilde{C}$ be the normalization of $C$, and $\tilde{\E}$ the pullback
of $\E$ to $\tilde{C}$. Given a logarithmic connection $\nabla$ on $\E$, we 
get a $D_C$-logarithmic connection $\tilde{\nabla}$ on $\tilde{\E}$, where $D_C$
is the divisor of points lying above the nodes of $C$.
\end{notn}

We want a complete
description of connections on $\tilde{\E}$ arising this way, and a
correspondence between these and connections on $\E$. We claim:

\begin{prop}\label{def-glue-basic} Logarithmic connections $\nabla$ on $\E$ 
are equivalent to connections
on $\tilde{\E}$ having simple poles at the points $P_1, Q_1,
\dots, P_\delta, Q_\delta$ lying above the nodes of $C$, and such that under 
the gluing
maps $G_i : \tilde{\E}|_{P_i} \rightarrow \tilde{\E}|_{Q_i}$ giving $\E$, for
each $i$ we have $\Res_{P_i}(\nabla) = - G_i^{-1} \circ \Res _{Q_i}(\nabla)
\circ G_i$. The properties of having trivial determinant and vanishing
$p$-curvature are preserved under this correspondence.
\end{prop}

\begin{proof} The main assertion follows easily from
\cite[Thm. 5.2.3]{co2} together with the remark \cite[p. 226]{co2} for 
nodal curves, which together state that sections of $\omega_C$ correspond
to sections of $\Omega^1_{\tilde{C}}(D_C)$ with residues at the pair of 
points above any given node adding to zero.

Since vanishing $p$-curvature can be verified on open sets, and the
normalization map is an isomorphism away from the nodes, it is clear that logarithmic 
connections with
vanishing $p$-curvature on $C$ will correspond to logarithmic connections 
with vanishing $p$-curvature on $\tilde{C}$. The same argument also works
for trivial determinant.
\end{proof}

We can in particular conclude:

\begin{cor}\label{def-conn-line} Let $\L$ be a line bundle on $C$. Then $\L$ 
can have a logarithmic
connection $\nabla$ with vanishing $p$-curvature only if $p | \deg
\tilde{\L}$.
\end{cor}

\begin{proof}Applying the previous proposition, if we pull back to
$\tilde{\nabla}$ on $\tilde{\L}$ we find that the residues of
$\tilde{\nabla}$ come in additive inverse pairs mod $p$. We obviously have
$p | \F^* (\tilde{\L}^{\tilde{\nabla}})$, and then by 
\cite[Cor. 2.11]{os10} we have that the determinant of the inclusion map
$\F^* (\tilde{\L}^{\tilde{\nabla}}) \hookrightarrow \tilde{\L}$
has total order equal to the sum of the residues mod $p$, which is zero, so
we conclude that $\deg \tilde{\L}$ must also vanish mod $p$, as asserted.
\end{proof}

We now restrict to the situation:

\begin{sit} Suppose that $\E$ has rank 2 and trivial determinant, 
and we have fixed an exact sequence 
$$0 \rightarrow \L \rightarrow \E \rightarrow \L^{-1} \rightarrow 0.$$ 
The same statements then hold for $\tilde{\E}$.
\end{sit}

We introduce some terminology in this situation:

\begin{defn} Given a logarithmic connection $\nabla$ on $\E$ 
(resp., a $D$-logarithmic connection $\nabla$ on $\tilde{\E}$), 
the {\bf Kodaira-Spencer map} associated to $\nabla$ and a
sub-line-bundle $\L$ (resp., $\tilde{\L}$) is the natural map
$\L \rightarrow \L^{-1} \otimes \omega_C$ (respectively, $\tilde{\L} 
\rightarrow \tilde{\L}^{-1} \otimes \Omega^1_{\tilde{C}}(D)$) obtained by 
composing the inclusion map, $\nabla$, and the quotient map. One verifies directly that this is a linear map. 

In the case that $\E$ (resp., $\tilde{\E}$) is unstable, we will refer to the Kodaira-Spencer map of $\nabla$ to mean the map associated to $\nabla$ and its destabilizing line bundle.
\end{defn}

Recall that by Lemma \ref{exp-destab-unique}, the destabilizing line bundle is unique, so the last part of the definition is justified. Note that with this terminology, Joshi and Xia's proof of \ref{exp-unstable} boils down to the statement that the Frobenius-pullback
of a Frobenius-unstable bundle necessarily has a connection such that
the Kodaira-Spencer map of the destabilizing line bundle is an isomorphism. 
It should perhaps therefore not be surprising that we will consider 
connections for which the Kodaira-Spencer is an isomorphism. We note:

\begin{lem}\label{def-arith-genus} Suppose that the arithetmic genus of $C$ 
(resp., the genus 
of $\tilde{C}$ plus $\frac{\deg D}{2}$) is greater than or equal to $3/2$;
that is to say, we are in the ``hyperbolic'' case. Then if the
Kodaira-Spencer map associated to $(\nabla, \L)$ is an isomorphism for any
$\nabla$, then $\L$ is a destabilizing line bundle for $\E$ (resp.,
$\tilde{\E}$), and is thus uniquely determined even independent of $\nabla$.
\end{lem}

\begin{proof} The Kodaira-Spencer isomorphism gives 
$\L^{\otimes 2} \cong \omega_C$ (resp., 
$\L^{\otimes 2} \cong \Omega^1_{\tilde{C}}(D)$), which from the hypotheses has positive degree.
\end{proof}

One can approach the issue of gluing connections from two perspectives:
either fixing the glued bundle $\E$ on $C$, and exploring which connections
on $\tilde{\E}$ will glue to yield connections on $\E$, or allowing the
gluing of $\E$ itself to vary as well. The author had originally intended to use the first approach, since we ultimately wish to classify the
connections on a particular unstable bundle on a nodal curve. 
However, the second approach, pursued by Mochizuki \cite[p. 118]{mo3}, offers a more transparent view of the more
general setting, and ultimately yields a cleaner argument even for our
specific application. As such, we now fix $\tilde{\E}$ on
$\tilde{C}$, but do not assume a fixed gluing $\E$ on $C$. That is to say:

\begin{sit}\label{def-glue} Fix $\tilde{\E}$ of rank 2 and trivial 
determinant, together with an exact sequence 
$$0 \rightarrow \tilde{\L} \rightarrow \tilde{\E} \rightarrow 
\tilde{\L}^{-1} \rightarrow 0.$$ 
\end{sit}

The main statement on gluing is:

\begin{prop}\label{def-glue-prop} In Situation \ref{def-glue},
let $\tilde{\nabla}$ be a $D_C$-logarithmic connection on $\tilde{C}$
with trivial determinant and vanishing $p$-curvature, such that the
Kodaira-Spencer map associated to $\tilde{\L}$ is an isomorphism. Further
suppose that the $e_1,e_2$ of \cite[Cor. 2.10]{os10} match one another (up 
to permutation) for pairs of points lying above given nodes of $C$. 
Then if one fixes a gluing $\L$ of $\tilde{\L}$ with $\L^{\otimes 2} \cong \omega_C$, there is a unique gluing of
$(\tilde{\E}, \tilde{\nabla})$ to a pair $(\E, \nabla)$ on $C$, such that
one obtains a sequence
$$0 \rightarrow \L \rightarrow \E \rightarrow \L^{-1} \rightarrow 0,$$ 
and the resulting $(\E, \nabla)$ will also have Kodaira-Spencer map an
isomorphism. If $C$ has arithmetic genus at least $2$, transport equivalence 
is preserved under this correspondence.
\end{prop}

\begin{proof} We first claim that the condition that the Kodaira-Spencer map
for $\tilde{\L}$ be an isomorphism implies that for any $P \in \{P_i,
Q_i\}$, $\tilde{\L}|_P$ is not contained in an eigenspace of $\Res_P
\tilde{\nabla}$, and that the eigenvalues are both non-zero. But due to
the triviality of the determinant, the sum of the eigenvalues is zero, so 
because the residue matrices are diagonalizable (see \cite[Cor. 2.11]{os10}),
the latter assertion is actually a consequence of the former. Now, 
considering the definition of the Kodaira-Spencer map $\tilde{\L} 
\rightarrow \tilde{\L}^{-1}\otimes \Omega^1_{\tilde{C}}(D_C)$, if we
restrict to $P$ we get a map which is clearly equal to zero if and only if
$\nabla(\tilde{\L}) |_P \subset \tilde{\L}\otimes \Omega^1_{\tilde{C}}|_P$, 
which is the case precisely
when $\tilde{\L}|_P$ is contained in an eigenspace of 
$\Res_P \tilde{\nabla}$, as desired.

Given this, for each pair $P_i, Q_i$, Proposition \ref{def-glue-basic} and 
our hypothesis on the matching eigenvalues of the residue matrices at $P_i,
Q_i$ imply that in order 
to glue the connection, it is necessary and sufficient to map eigenspaces of 
opposing sign to each other. To glue $\tilde{\L}$, we also map its image
at $P_i$ to its image at $Q_i$. We thus see that the two eigenspaces of 
$\Res_{P_i} \tilde{\nabla}$ and
$\Res_{Q_i} \tilde{\nabla}$ and the images of $\tilde{\L}$ form a set of
three one-dimensional subspaces which must be matched under $G_i$, and it is
easy to see that this determines $G_i$ up to scaling. But finally, scaling of
$G_i$ is equivalent to scaling the induced gluing map on $\tilde{\L}$, which
is precisely what determines the isomorphism class of the glued $\L$; thus,
$\L$ may be specified arbitrarily, and given a choice of $\L$, the $G_i$ and
hence the pair $(\E, \nabla)$ are uniquely determined, as desired. Lastly,
we observe that since the Kodaira-Spencer map gives an isomorphism $\L
\otimes (\E/\L)^{-1} \cong \omega_C$, the hypothesis that $\L^{\otimes 2} \cong \omega_C$ is equivalent to the condition that the glued $\E$ have trivial determinant.

Considering transport, it is trivial that if two connections on 
$\E$ are transport-equivalent, then their pullbacks to $\tilde{\E}$ are, 
and for the converse, the uniqueness of the gluing makes it clear that if
two connections $\tilde{\nabla}$ and $\tilde{\nabla'}$ on $\tilde{\E}$ are 
transport-equivalent under an automorphism $\vp$ of $\tilde{\E}$, then 
$\vp$ naturally gives an isomorphism of the two gluings $\E$ and $\E'$
which takes $\nabla$ to $\nabla'$. Finally, the hypothesis that the 
arithmetic genus of $C$ is at least $2$ implies that $\L$ and $\tilde{\L}$ 
are uniquely determined as the destabilizing sub-bundles of $\E$ and 
$\tilde{\E}$, so there is no concern that they might change under 
transport.
\end{proof}

Putting together the previous propositions, we finally conclude:

\begin{cor}\label{def-glue-main} Let $\tilde{\E}$ be a vector bundle on
$\tilde{C}$ of rank 2, with the arithmetic genus of $C$ being at least 2,
and suppose there exists an exact sequence
$$0 \rightarrow \tilde{\L} \rightarrow \tilde{\E} \rightarrow
\tilde{\L}^{-1} \rightarrow 0.$$
Fix a gluing of $\tilde{\L}$ to a line bundle $\L$ on $C$ satisfying $\L^{2} \cong \omega_C$. Then there
exists a bijective equivalence between transport-equivalence classes of 
$D_C$-logarithmic connections $\tilde{\nabla}$ on $\tilde{\E}$ with trivial
determinant and vanishing $p$-curvature, the eigenvalues of the residues of
$\tilde{\nabla}$ matching at the pairs of points above each node, and having 
the Kodaira-Spencer map an isomorphism on one side, and on the other
side, pairs $(\E, \nabla)$ of gluings of $\E$ preserving an exact sequence
$$0 \rightarrow \L \rightarrow \E \rightarrow \L^{-1} \rightarrow 0,$$
together with logarithmic connections $\nabla$ on $\E$ with vanishing 
$p$-curvature and trivial determinant and having the Kodaira-Spencer map an
isomorphism, up to isomorphism and transport equivalence.

Further, this correspondence holds for first-order infinitesmal deformations.
\end{cor}

\begin{proof} We can immediately conclude the statement over a field from
our previous propositions. For first-order deformations, the same arguments
will go through, with the aid of the following facts: first and most
substantively, it follows from \cite[Cor. 3.6]{os10} that 
the residue matrices on $\tilde{C}$ will still be diagonalizable 
over $k[\epsilon]/\epsilon^2$, with the eigenvalues $e_i$ the same as for
the connection being deformed. Next, since we are simply taking a base
change of our original situation over $k$, the general gluing description 
given by Proposition \ref{def-glue-basic} still holds for formal reasons. 
Finally, one can easily verify that even over an arbitrary ring, it is still the case that an automorphism of a rank two free module is 
determined uniquely by sending any three pairwise independent lines to any 
other three. We therefore conclude the desired statement for first-order
deformations as well. 
\end{proof}

\section{Deforming to a Smooth Curve}\label{s-def-deform}

The ultimate goal of this section is to prove that the connections we are
interested in can always be smoothed from a general irreducible rational
nodal curve, which together with the finiteness result of Section
\ref{s-exp-det} and the main results of \cite{os7}, \cite{os10}, will allow us to finish the proof of the characteristic-independent portion of Theorem \ref{exp-main}. We begin with some general observations on when the space of connections with vanishing $p$-curvature is smooth over a given deformation of the curve and vector bundle. We then make a key dimension computation using the 
techniques of \cite{os10} and of the previous section, once again following arguments of Mochizuki \cite[Cor. II.2.5, p. 150]{mo3} rather than the original approach of the author, for the sake simplicity and generality.

\begin{sit}\label{def-sit}
We suppose that $C_0$ is an irreducible, rational proper curve with two 
nodes, $\tilde{C_0} \cong \P^1$ its normalization, with $P_1, Q_1, P_2, Q_2$ 
being the points lying above the two nodes. We let $\E_0$ be the vector bundle
described by Situation \ref{exp-specific-e}, and $\nabla_0$ a logarithmic 
connection on $\E_0$ with trivial determinant and vanishing $p$-curvature.
\end{sit}

By Proposition \ref{def-background}, $p$-curvature 
gives an algebraic morphism
$\psi_p : H^0(\cEnd^0(\E_0) \otimes \omega_{C_0})
\rightarrow H^0(\cEnd^0(\E_0) \otimes F^* \omega_{C_0^{(p)}})$
such that for $\vp \in H^0(\cEnd^0(\E_0) \otimes \omega_{C_0})$,
$\psi_p (\nabla_0 + \vp)$ in fact lies in 
$H^0(\cEnd^0(\E_0) \otimes F^* \omega_{C_0^{(p)}})^{(\nabla_0+\vp)^{\ind}}$. Now, 
we first claim:

\begin{lem}\label{def-dpsi} If $\nabla_0$ has vanishing $p$-curvature, the differential of 
$\psi_p$ at $0$ gives a linear map
$$d\psi_p: H^0(\cEnd^0(\E_0) \otimes \omega_{C_0}) \rightarrow 
H^0(\cEnd^0(\E_0) \otimes F^* \omega_{C_0^{(p)}})^{\nabla_0^{\ind}}.$$
\end{lem}

\begin{proof} We simply consider the induced map on first-order deformations
of $\nabla_0$. Denoting for the moment by $C_1$, $\E_1$ the base
change of $C_0$, $\E_0$ to $k[\epsilon]/(\epsilon^2)$, suppose that 
$\vp \in \epsilon H^0(\cEnd^0(\E_1) \otimes \omega_{C_1}) \cong
H^0(\cEnd^0(\E_0) \otimes \omega_{C_0})$, 
and consider $\nabla_0 + \vp$. Since $\nabla_0$ has vanishing
$p$-curvature, the image under $\psi_p$ is in 
$\epsilon H^0(\cEnd^0(\E_1) \otimes F^* \omega_{C_1^{(p)}})
^{(\nabla_0+ \epsilon \vp)^{\ind}}$, which is naturally isomorphic to 
$H^0(\cEnd^0(\E_0) \otimes F^* \omega_{C_0^{(p)}})^{\nabla_0^{\ind}}$, giving 
the desired result.
\end{proof}

Our main assertion is:

\begin{prop} If the map $d\psi_p$ of the previous lemma is surjective, then
given a deformation $C$ of $C_0$ and $\E$ of $\E_0$ on $C$, such that
the functor of connections on $\E$ with trivial determinant is formally 
smooth at $\nabla_0$, then the functor of connections on $\E$ with trivial 
determinant and vanishing $p$-curvature is formally smooth at $\nabla_0$. 
\end{prop}

\begin{proof} By hypothesis, there is no obstruction to deforming 
$\nabla_0$ as a connection with trivial determinant. 
Following \cite[Def. 1.2,
Rem. 2.3]{sc2}, we say that a map $B \twoheadrightarrow A$ of local Artin 
rings over the base ring of our
deformation and having residue field $k$ is a {\bf small extension} if 
the kernel is a principal ideal $(\epsilon)$ with $(\epsilon) \m_B = 0$; 
it follows then that $\epsilon B \subset B$ is isomorphic to $k$. To 
verify (formal) smoothness, by virtue of \cite[Prop. 17.14.2]{ega44} it is 
easily checked inductively that it is 
enough to check on small extensions. We show therefore that for such a
small extension, when $d\psi_p$ is surjective there is no obstruction to 
lifting a deformation of $\nabla_0$ over $A$ to a deformation over $B$,
even with the addition of the vanishing $p$-curvature hypothesis.
Let $C_B, \E_B$ be the given deformations over $B$ of $C_0, \E_0$ 
respectively, with $C_A, \E_A$ the induced deformations over $A$, and 
suppose that $\nabla_B$ is a connection on $\E_B$ such that $\nabla_A$ 
has vanishing $p$-curvature. The main point is that it is straightforward
to check that the hypothesis that $\epsilon B \cong k$ implies that
$\epsilon H^0(\cEnd^0(\E_B) \otimes \omega_{C_B}) \cong 
H^0(\cEnd^0(\E_0) \otimes \omega_{C_0})$, and for any
$\vp \in \epsilon H^0(\cEnd^0(\E_B) \otimes \omega_{C_B})$, we have
$\epsilon H^0(\cEnd^0(\E_B) \otimes F^* \omega_{C_B^{(p)}})
^{(\nabla_B+ \vp)^{\ind}} \cong H^0(\cEnd^0(\E_0) \otimes 
F^* \omega_{C_0^{(p)}})^{\nabla_0^{\ind}}$.
We want to show that for some choice of 
$\vp \in \epsilon H^0(\cEnd^0(\E_B) \otimes \omega_{C_B})$,
$\nabla_B + \vp$ has vanishing $p$-curvature. But as before, since 
$\nabla_A$ has vanishing $p$-curvature, the image under $\psi_p$ of 
$\nabla_B + \vp$ is in $\epsilon H^0(\cEnd^0(\E_B) \otimes F^* \omega_{C_B^{(p)}})
^{(\nabla_B+ \vp)^{\ind}}$, and under the above isomorphisms, the 
induced map is equal to $d\psi_p + \frac{1}{\epsilon} \psi_p(\nabla_B)$, 
where $\frac{1}{\epsilon}$ is simply shorthand for the isomorphism 
$\epsilon H^0(\cEnd^0(\E_B) \otimes F^* \omega_{C_B^{(p)}}) 
^{(\nabla_B+ \vp)^{\ind}}\risom H^0(\cEnd^0(\E_0) \otimes F^* 
\omega_{C_0^{(p)}})^{\nabla_0^{\ind}}$. Hence if $d \psi_p$ is
surjective, we can choose $\vp$ so that $\nabla_B + \vp$ has
vanishing $p$-curvature, as desired. 
\end{proof}

We observe that in our situation, the normalization $\tilde{\E_0}$ of $\E_0$ 
is isomorphic to $\O(1) \oplus \O(-1)$: we certainly have $\tilde{\L}\cong
\O(1)$, so by Lemma \ref{exp-destab-unique},
$\O(1)$ is the maximal line bundle in $\tilde{\E_0}$, and then the desired
splitting follows from \cite[Proof of Thm. 1.3.1]{h-l}.
Also, by Proposition \ref{def-glue-basic} $\tilde{\nabla_0}$ is a
$D_{C_0}$-logarithmic connection on $\tilde{\E_0}$ with trivial determinant 
and vanishing $p$-curvature. For the sake of cleanness and generality, we
use Mochizuki's arguments \cite[Cor. II.2.5, p. 150]{mo3} to prove the
following.

\begin{prop}\label{def-aux} If $\nabla_0$ has a non-zero Kodaira-Spencer
map, then the space of sections of $\cEnd^0(\E_0) \otimes F^*
\omega_{C_0^{(p)}}$
horizontal with respect to the connection $\nabla^{\ind}_0$ induced by 
$\nabla_0$ on $\E_0$ and $\Nc$ on $F^* \omega_{C_0^{(p)}}$ has 
dimension $3$.
\end{prop}

\begin{proof} The proof proceeds in two parts: we show that $H^1(C_0,
(\cEnd^0(\E_0) \otimes F^* \omega_{C_0^{(p)}})^{\nabla_0^{\ind}})=0$, and then
compute the Euler characteristic. Both computations require formal local
computations, so we begin by setting out the situation formally locally at a
node of $C_0$. First, note that although taking kernels and tensor products of
connections do not commute in general, there is no problem when one
connection is obtained as the canonical connection of a Frobenius pullback,
so we have $(\cEnd^0(\E_0) \otimes \F^* \omega_{C_0^{(p)}})^{\nabla_0^{\ind}} =
\cEnd^0(\E_0)^{\nabla_0^{\cEnd}} \otimes \omega_{C_0^{(p)}}$. Formally locally
at the node, $C_0$ is isomorphic to $k[[x,y]]/(x,y)$; moreover, we claim
that if we choose $x,y$ correctly, we can trivialize $\E_0$ so that
$\nabla_0^{\cEnd}$ has connection matrix $\begin{bmatrix}e (\frac{dx}{x} - \frac{dy}{y}) & 0 & 0 \\
0 & 0 & 0 \\
0 & 0 & e (\frac{dy}{y} - \frac{dx}{x}) \end{bmatrix}$ for some $e$ with $0<e<p$. Indeed, this follows from Proposition \ref{def-glue-basic} together with the formal 
local diagonalizability result of 
\cite[Cor. 2.10]{os10} applied to $\tilde{C}_0$, noting that if the residue of
$\nabla_0$ has eigenvalues $e', -e'$, then the residue of $\nabla_0^{\cEnd}$ has eigenvalues $2e', 0, -2e'$. By the same token, the pullback to the normalization gives connection matrices
$\begin{bmatrix}e \frac{dx}{x} & 0 & 0 \\
0 & 0 & 0 \\
0 & 0 & -e \frac{dx}{x}\end{bmatrix}$ and $\begin{bmatrix}-e \frac{dy}{y} & 0 & 0 \\
0 & 0 & 0 \\
0 & 0 & e \frac{dy}{y}\end{bmatrix}$. Finally, we note that the kernel of
the connection on $C_0$ is given over $\O_{C_0^{(p)}}$ by $(x^{p-e}, y^e) \oplus
(1) \oplus (x^e, y^{p-e})$, and by $(x^{p-e}) \oplus (1) \oplus (x^{e})$ and
$(y^{e}) \oplus (1) \oplus (y^{p-e})$ on the normalization. The formal local
calculations of the following paragraphs are justified by the following
facts: given a sheaf map, surjectivity, and more generally factoring through
a given subsheaf, may be checked after completion; completion commutes with
pullback, with taking kernels of connections in characteristic $p$, and with
modding out by torsion over a DVR; finally, completion is well-behaved with
respect to pushforward under the normalization map by the theorem on formal
functions. 

Now, to check that $H^1$ vanishes, by Grothendieck duality on $C_0^{(p)}$ it
suffices to check that $\Hom(\cEnd^0(\E_0)^{\nabla_0^{\cEnd}} \otimes
\omega_{C_0^{(p)}}, \omega_{C_0^{(p)}}) = \Hom(\cEnd^0(\E_0)^{\nabla_0^{\cEnd}},
\O_{C_0^{(p)}})=0$. Although a section of the latter need not come from a map 
$\cEnd^0(\E_0)\rightarrow \O_{C_0}$ which is horizontal with respect to 
$\nabla_0^{\cEnd}$, we claim that it does after normalization. We have a
natural map $\cHom(\cEnd^0(\E_0)^{\nabla_0^{\cEnd}},
\O_{C_0^{(p)}})|_{\tilde{C}_0} \rightarrow
\cHom(\cEnd^0(\E_0)|_{\tilde{C}_0}^{\nabla_0^{\cEnd}}, \O_{\tilde{C}_0^{(p)}})$,
and a natural inclusion $\cHom(\cEnd^0(\E_0)|_{\tilde{C}_0},
\O_{\tilde{C}_0})^{\nabla_0^{\cEnd}} \hookrightarrow
\cHom(\cEnd^0(\E_0)|_{\tilde{C}_0}^{\nabla_0^{\cEnd}},
\O_{\tilde{C}_0^{(p)}})$. These are both isomorphisms away from the points 
above the nodes, for trivial reasons in the first case, and because of
Theorem \ref{exp-cartier-vect} for the second. We want to show
that the first map factors through the second. Examining the formal local
situation at a node, we first note that if $e_1, e_2>0$, any map from
$(x^{e_1}, y^{e_2})$ to $\O_{C_0^{(p)}}$ necessarily vanishes, and more specifically, sends
$x^{e_1}$ and $y^{e_2}$ to positive ($p$th) powers of $x$ and $y$
respectively. It is thus clear that give a map $\cEnd^0(\E_0)^{\nabla_0^{\cEnd}}
\rightarrow \O_{C_0^{(p)}}$, after normalization we can divide through to get
a map formally locally $\cEnd^0(\E_0)|_{\tilde{C}_0} \rightarrow
\O_{\tilde{C}_0}$
which commutes with the induced connection, completing the proof of the
claim. Next, we claim that such a map must be $0$. Indeed, if we consider
the line sub-bundle $\L^0 \subset \cEnd^0(\E_0)$ which sends $\L \subset
\E_0$
to $0$, we see that it is isomorphic to $\L^{\otimes 2}$, and is not horizontal for
$\nabla_0^{\cEnd}$, since $\L$ is not horizontal for $\nabla_0$, and we have the
same situation after normalization. But having such a destabilizing line
sub-bundle precludes the existence of a horizontal morphism
$\cEnd^0(\E_0)|_{\tilde{C}_0} \rightarrow \O_{\tilde{C}_0}$ by Proposition \ref{exp-destab-unique}, so we conclude the desired vanishing statement.

Thus, it remains to compute the Euler characteristic of
$\cEnd^0(\E_0)^{\nabla_0^{\cEnd}}$. Since we only have two non-zero eigenvalues at each $P_i$ or $Q_i$, it follows from \cite[Cor. 2.11]{os10} that the cokernel of
$F^*((\cEnd^0(\E_0)|_{\tilde{C}_0})^{\nabla_0^{\cEnd}}) \rightarrow
\cEnd^0(\E_0)|_{\tilde{C}_0}$ is supported at the $P_i,Q_i$, with length $p$ at
each point. Since $\deg(\cEnd^0(\E_0)|_{\tilde{C}_0})=0$, we find that
$\deg(F^*((\cEnd^0(\E_0)|_{\tilde{C}_0})^{\nabla_0^{\cEnd}}))=-4p$. Next, we claim
that $(\cEnd^0(\E_0)|_{\tilde{C}_0})^{\nabla_0^{\cEnd}}$ is isomorphic to the
quotient of $(\cEnd^0(\E_0)^{\nabla_0^{\cEnd}})|_{\tilde{C}_0}$ by its torsion,
which we denote by $\F$; indeed, we clearly have a morphism from the latter
to the former, which is an isomorphism away from the points above the nodes,
hence gives an injection since we modded out by torsion. Surjectivity above
the nodes is then checked formally locally from our above description, so we
have $\deg(\F) = -4$, and $\chi(\F)=-1$. Finally, if $\nu$ denotes the
normalization map, we claim that the natural injection
$\cEnd^0(\E_0)^{\nabla_0^{\cEnd}} \hookrightarrow \nu_* \F$ has cokernel of
length $1$ at each node; again, this is checked formally locally, noting
that the cokernel will arise only from the summand at each node on which the
connection vanishes. We conclude therefore that $\chi
(\cEnd^0(\E_0)^{\nabla_0^{\cEnd}}) = -3$, so $H^0
(\cEnd^0(\E_0)^{\nabla_0^{\cEnd}}\otimes \omega_{C_0^{(p)}}) = \chi
(\cEnd^0(\E_0)^{\nabla_0^{\cEnd}}\otimes \omega_{C_0^{(p)}}) = 3$, completing the
proof of the proposition. 
\end{proof}

Finally, we put these results together in our specific situation:

\begin{thm}\label{def-deform} Let $C_0$ be a nodal rational curve of genus $2$, and 
$\E_0$ as in Situation \ref{exp-specific-e}. Let $\nabla_0$ have vanishing 
$p$-curvature and trivial determinant, and suppose that $\nabla_0$ has
no deformations preserving the $p$-curvature and not arising from transport. 
Then the map $d\psi_p$ of Lemma
\ref{def-dpsi} is surjective; in particular, given any deformation $C$ of
$C_0$, if $\E$ is the corresponding deformation of $\E_0$, then
the space of connections with trivial determinant and vanishing $p$-curvature 
on $\E$ is formally smooth at $\nabla_0$.
\end{thm}

\begin{proof} The main point is that by Remark \ref{exp-gen}, the space
of transport-equivalence classes of connections with trivial determinant 
on $\E_0$ or $\E$ is explicitly parametrized by $\A^3$ over the appropriate
base. In particular, deformations of $\nabla_0$ as a connection with trivial
determinant are unobstructed, and it also follows that
the space of first-order deformations of $\nabla_0$ with trivial
determinant, modulo those arising from transport, is three-dimensional.
By Proposition \ref{def-aux}, the image space of $d\psi_p$ is
three-dimensional. We therefore get surjectivity precisely when transport 
accounts for the entire kernel, which is to say, when there are no 
deformations of $\nabla_0$ having vanishing $p$-curvature and trivial 
determinant other than those obtained by transport. We can thus apply the 
previous proposition to conclude smoothness.
\end{proof}

It is now a matter of some simple combinatorics to complete the proof of the characteristic-independent portion of Theorem \ref{exp-main}.

\begin{proof}[Proof of Theorem \ref{exp-main}, $p>2$ case]
By the results of Section \ref{s-exp-background} it suffices to show
that, for the particular $\E$ of Situation \ref{exp-specific-e},
there are precisely $\frac{1}{24}p(p^2-1)$ transport-equivalence classes of
connections with trivial determinant and vanishing $p$-curvature on $\E$,
and that none of these have any non-trivial deformations. We will show that
this statement holds in the situation that $C$ is a general rational nodal
curve, and then conclude the same result must hold for a general smooth
curve.

We observe that even in the situation of a nodal curve, there is a unique
extension $\E$ of $\L^{-1}$ by $\L$; indeed, the proof of Proposition
\ref{exp-unstable-unique} goes through with $\omega_C$ in place of
$\Omega^1_C$. We also note that by Corollary
\ref{def-conn-line}, the argument of Proposition \ref{exp-unstable} still 
shows that any
connection must have its Kodaira-Spencer map be an isomorphism. It then 
follows from Corollary \ref{def-glue-main} that it suffices
to prove the same result for $D$-logarithmic connections on $\O(1) \oplus
\O(-1)$ on $\P^1$ satisfying the hypotheses of 
\cite[Sit. 2.12]{os10} and having the Kodaira-Spencer map
an isomorphism, where $D$ is made up of four 
general points on $\P^1$, and the eigenvalues of the residues at the points
match in the appropriate pairs. We note that by degree considerations, the
Kodaira-Spencer map in this case is always either zero or an isomorphism,
so if we fix eigenvalues $\alpha_i$ for each pair $(P_i, Q_i)$, by 
\cite[Thm. 1.1]{os10} we find that we are
looking for separable rational functions on $\P^1$ of degree $2p-1 
- 2 \sum \alpha_i$, and ramified to order at least 
$p- 2\alpha_i$ at $P_i$ and $Q_i$ (note that the coefficient 
doubling for the degree is due to our use of a single, matching $\alpha_i$ 
for both
$P_i$ and $Q_i$). We could use the second formula of \cite[Cor. 8.1]{os7} 
to compute the answer directly, but the first formula yields
a more elegant solution. In either case, we are already given the lack of
non-trivial deformations, so it suffices to show that the number of maps is
correct. The formula gives that for each $(\alpha_1, \alpha_2)$ there are
$$\min\{\{p-2 \alpha_i\}_i, \{p-2 \alpha_{3-i}\}_i, \{2\alpha_i\}_i,
\{2\alpha_{3-i}\}_i\}$$
such maps, which reduces to
$$\min\{\{p-2 \alpha_i\}_i, \{2 \alpha_i\}_i\}.$$
Rather than summing up over all $\alpha_i$, as we would with the second
formula, we note that the number of maps will also be given by:
$$\sum_{1 \leq j \leq (p-1)/2} \#\{(\alpha_1, \alpha_2): j \leq 2 \alpha_i, j
\leq p -2 \alpha_i\}$$ 
which then reduces to
$$ \sum_{1 \leq j \leq (p-1)/2} \!\!\!\!\!\!\!\! (\frac{p+1}{2}-j)^2 =
\!\!\!\!\!\!\!\! \sum_{1 \leq j \leq (p-1)/2} \!\!\!\!\!\!\!\! j^2 =
\sum_{1 \leq j \leq (p-1)/2} \!\!\!\!\!\!\!\! (2
\binom{j}{2} + j)$$ 
$$= 2 \binom{(p+1)/2}{3} + \frac{p+1}{2}\frac{p-1}{4} =
\frac{1}{24}(p+1)((p-1)(p-3)+ 3(p-1)) = \frac{(p+1)(p-1)p}{24},$$
giving the desired result for a general nodal curve.

We can now apply Theorem \ref{def-deform} to conclude that since none of our
connections on the general nodal curve have non-trivial deformations, the 
space of
connections with trivial determinant and vanishing $p$-curvature on our
chosen bundle over our parameter space of genus $2$ curves is formally
smooth at each connection on the general nodal curve. Furthermore, by
Corollary \ref{exp-finite} (in light of Remark \ref{exp-gen}), this space of
connections is finite, so we conclude that it is finite \'etale at the general
nodal curve, and finite everywhere, which then implies (i) for a general
smooth curve, as desired.
\end{proof}

\bibliographystyle{hamsplain}
\bibliography{hgen}

\newcommand{\noopsort}[1]{} \newcommand{\printfirst}[2]{#1}
  \newcommand{\singleletter}[1]{#1} \newcommand{\switchargs}[2]{#2#1}
\providecommand{\bysame}{\leavevmode\hbox to3em{\hrulefill}\thinspace}
\begin{thebibliography}{10}

\bibitem{co2}
Brian Conrad, \emph{Grothendieck duality and base change}, Lecture Notes in
  Mathematics, no. 1750, Springer-Verlag, 2000.

\bibitem{gi2}
D.~Gieseker, \emph{Stable vector bundles and the {F}robenius morphism}, Ann.
  Sci. Ecole Norm. Sup. (4) \textbf{6} (1973), 95--101.

\bibitem{ha1}
Robin Hartshorne, \emph{Algebraic geometry}, Springer-Verlag, 1977.

\bibitem{h-l}
Daniel Huybrechts and Manfred Lehn, \emph{The geometry of moduli spaces of
  sheaves}, Max-Planck-Institut fur Mathematik, 1997.

\bibitem{j-r-x-y}
Kirti Joshi, S.~Ramanan, Eugene~Z. Xia, and Jiu-Kang Yu, \emph{On vector
  bundles destabilized by {F}robenius pull-back}, preprint.

\bibitem{j-x}
Kirti Joshi and Eugene~Z. Xia, \emph{Moduli of vector bundles on curves in
  positive characteristic}, Compositio Math. \textbf{122} (2000), no.~3,
  315--321.

\bibitem{ka1}
Nicholas~M. Katz, \emph{Nilpotent connections and the monodromy theorem:
  Applications of a result of {Turrittin}}, Inst. Hautes Etudes Sci. Publ.
  Math. \textbf{39} (1970), 175--232.

\bibitem{l-p3}
Herbert Lange and Christian Pauly, \emph{On {F}robenius-destabilized rank-2
  vector bundles over curves}, \mbox{arXiv:math.AG/0309456}.

\bibitem{l-p}
Yves Laszlo and Christian Pauly, \emph{The action of the {F}robenius maps on
  rank 2 vector bundles in characteristic 2}, Journal of Algebraic Geometry
  \textbf{11} (2002), no.~2, 129--143.

\bibitem{mo1}
Shinichi Mochizuki, \emph{A theory of ordinary $p$-adic curves}, Publ. RIMS
  \textbf{32} (1996), no.~6, 957--1151.

\bibitem{mo3}
\bysame, \emph{Foundations of $p$-adic {T}eichm\"uller theory}, American
  Mathematical Society, 1999.

\bibitem{os9}
B.~Osserman, \emph{The generalized {V}erschiebung map for curves of genus 2},
  preprint.

\bibitem{os10}
\bysame, \emph{Logarithmic connections with vanishing $p$-curvature},
  \mbox{arXiv:math.AG/0409145}.

\bibitem{os6}
\bysame, \emph{Mochizuki's crys-stable bundles: A lexicon and applications},
  preprint.

\bibitem{os7}
\bysame, \emph{Rational functions with given ramification in characteristic
  $p$}, \mbox{arXiv:math.AG/0407445}.

\bibitem{ra2}
M.~Raynaud, \emph{Sections des fibr\'es vectoriels sur une courbe}, Bull. Soc.
  Math. France \textbf{110} (1982), no.~1, 103--125.

\bibitem{sc2}
Michael Schlessinger, \emph{Functors of {A}rtin rings}, Transactions of the AMS
  \textbf{130} (1968), 208--222.

\bibitem{ega44}
A.~Grothendieck with J.~Dieudonn\'e, \emph{{\'E}l\'ements de g\'eom\'etrie
  alg\'ebrique: {IV.} \'{E}tude locale des sch\'emas et des morphismes de
  sch\'emas, quatri\'eme partie}, vol.~32, Publications math\'ematiques de
  l'I.H.\'E.S., no.~2, Institut des Hautes \'Etudes Scientifiques, 1967.

\end{thebibliography}
\end{document}